\documentclass[a4paper]{amsart}

\PassOptionsToPackage{usenames,dvipsnames}{xcolor}
\usepackage{multirow,mathabx}
\usepackage[english]{babel}
\usepackage{fancyhdr}

\usepackage[utf8]{inputenc}

\usepackage{mathtools}
\usepackage{amsthm,amssymb,latexsym,url,cancel}
\usepackage{pdfpages}

\usepackage{tcolorbox}
\tcbuselibrary{theorems}

\usepackage{tikz}
\usepackage{tikz-qtree}
\usepackage{tkz-euclide}
\usepackage{pgfplots}
\usetkzobj{all}
\usetikzlibrary{intersections}
\usetikzlibrary{calc}
\usetikzlibrary{decorations.markings}
\usetikzlibrary{backgrounds,shapes}

\usepackage[all]{xy}
\usepackage{subfig}
\usepackage{enumerate}
\makeatletter
\let\@@enum@org\@@enum@
\def\@@enum@[#1]{\@@enum@org[\normalfont #1]}
\makeatother

\def\co{\colon\thinspace}

\DeclarePairedDelimiter{\ceil}{\lceil}{\rceil}

\DeclarePairedDelimiter{\brac}{[}{]}
\DeclarePairedDelimiter{\form}{\langle}{\rangle}
\DeclarePairedDelimiter{\fform}{\langle\!\langle}{\rangle\!\rangle}
\DeclarePairedDelimiter{\eval}{}{\rvert}

\DeclarePairedDelimiter{\abs}{\lvert}{\rvert}
\DeclarePairedDelimiter{\abss}{\lVert}{\rVert}
\newcommand\der[2]{#1^{(#2)}}
\DeclarePairedDelimiter{\mup}{\mu(}{)}
\newcommand\CR{\operatorname{CritReg}}

\newcommand\Inv{\operatorname{Inv}}

\newcommand\repart[1]{\mathrm{Re}\left[ #1 \right]}

\newcommand\ba{\begin{align*}}
\newcommand\ea{\end{align*}}
\newcommand\be{\begin{enumerate}}
\newcommand\ee{\end{enumerate}}
\newcommand\bp{\begin{proof}}
\newcommand\ep{\end{proof}}
\newcommand\bpp{\begin{prop}}
\newcommand\epp{\end{prop}}
\newcommand\bpb{\begin{prob}}
\newcommand\epb{\end{prob}}
\newcommand\bd{\begin{defn}}
\newcommand\ed{\end{defn}}
\newcommand\bh{\begin{hint}}
\newcommand\eh{\end{hint}}

\newcommand\bC{\mathbb{C}}

\newcommand\bN{\mathbb{N}}
\newcommand\bR{\mathbb{R}}

\newcommand\bZ{\mathbb{Z}}

\newcommand\BB{\mathcal{B}}
\newcommand\CC{\mathcal{C}}
\newcommand\DD{\mathcal{D}}
\newcommand\FF{\mathcal{F}}

\newcommand\GG{\mathcal{G}}
\newcommand\VV{\mathcal{V}}
\newcommand\UU{\mathcal{U}}

\newcommand\restrict{\!\!\restriction}

\newcommand\supp{\operatorname{supp}}
\newcommand\supt{\operatorname{supt}}
\newcommand\Id{\operatorname{Id}}

\newcommand\diam{\operatorname{diam}}

\newcommand\gam{\Gamma}
\newcommand\Mod{\operatorname{Mod}}

\DeclareMathOperator\Homeo{Homeo}

\newcommand\sse{\subseteq}

\DeclareMathOperator\hull{hull}
\DeclareMathOperator\Fix{Fix}

\DeclareMathOperator\Diff{Diff}

\DeclareMathOperator\PL{PL}

\DeclareMathOperator\dom{dom}

\renewcommand{\MR}[1]
{\href{http://www.ams.org/mathscinet-getitem?mr=#1}{MR#1}}

\usepackage{mathrsfs}

\usepackage{geometry}

\def\thetitle{Algebraic Structure of Diffeomorphism Groups of One--Manifolds}
\def\theshorttitle{Diffeomorphism Groups of $1$--Manifolds}
\usepackage{hyperref}
\hypersetup{
    pdftitle=   \thetitle,
   pdfauthor=  {S. Kim, T. Koberda, J. Chang}
}

\newtheorem{thm}{Theorem}[section]
\newtheorem{lem}[thm]{Lemma}

\newtheorem{cor}[thm]{Corollary}
\newtheorem{prop}[thm]{Proposition}

\newtheorem{que}[thm]{Question}
\newtheorem*{claim*}{Claim}

\theoremstyle{remark}
\newtheorem{exmp}[thm]{Example}

\newtheorem{rem}[thm]{Remark}
\newtheorem*{hint}{Hint}

\newtheorem{notation}[thm]{Notation}
\newtheorem*{soln}{Solution}

\def\bsl{\begin{soln}}
\def\esl{\end{soln}}

\theoremstyle{definition}
\newtheorem{defn}[thm]{Definition}
\newtheorem{prob}{Problem}[section]

\begin{document}

\title[\theshorttitle]\thetitle
\date{\today}

\author[J. Chang]{Jaewon Chang}
\address{School of Mathematics, Korea Institute for Advanced Study, Seoul, Korea}
\email{jjwjjw9595@gmail.com}
\author[S. Kim]{Sang-hyun Kim}
\address{School of Mathematics, Korea Institute for Advanced Study, Seoul, Korea}
\email{skim.math@gmail.com}
\urladdr{http://cayley.kr}

\author[T. Koberda]{Thomas Koberda}
\address{Department of Mathematics, University of Virginia, Charlottesville, VA 22904-4137, USA}
\email{thomas.koberda@gmail.com}
\urladdr{http://faculty.virginia.edu/Koberda}

\maketitle

\setcounter{tocdepth}{1}\tableofcontents


\section*{Introduction}
It is a celebrated result of Mather that the group of $C^k$--diffeomorphisms of an $n$--manifold is simple~\cite{Mather1,Mather2},
provided that a mild isotopy condition is satisfied, with the possible exception of $k=n+1$.
The purpose of this article is mostly expository, and in it we give a detailed account of Mather's proof in the case when $n=1$. 

The content of this article can be divided broadly into two parts.
The first part is to give a complete and self-contained account of the proof in the case of dimension $1$.
Mather's original articles assume the familiarity with some analytic and topological methods such as  concave moduli,
fragmentation and the Epstein--Ling Theorem. In the present manuscript, we concentrate on
 conveying the key ideas of his proof, including relevant background.
To achieve this first goal more efficiently, we focus on the one--dimensional case and avoid the notational confusion
 introduced by the case of general smooth manifolds.

The reader will find that at the technical heart of the proof lies in the
construction of a norm-reduction operator in Section~\ref{s:norm}.
In detailing this construction,
we will refine Mather's original result and establish the simplicity of $C^{k,\alpha}$ diffeomorphism groups of $\mathbb{R}$,
where $\alpha$ is a ``tame'' concave modulus. 
We will also indicate how the same ideas apply to general $n>1$ manifolds.

The second part is to survey some recent results of the first two authors
regarding isomorphism types of finitely generated subgroups of one-dimensional diffeomorphism groups. 
In particular, we deduce from our refined versions of Mather's results that homomorphisms
between one--dimensional diffeomorphism groups with relatively weak hypotheses on the relevant moduli of continuity
are very restricted, as observed in~\cite{KK2017crit}. 

\section{Main objectives}
A \emph{concave modulus (of continuity)} is a homeomorphism \[\alpha\co [0,\infty)\to[0,\infty)\] that is concave as a function.
For instance, the map $\omega_s(x):=x^s$ is a concave modulus for each $s\in(0,1]$. The moduli $\{\omega_s\}_{s\in(0,1]}$ are sometimes
called the \emph{H\"older} moduli.

Let $f\co X\to Y$ be a map between metric spaces. We set \[ [f]_\alpha=\sup_{x\ne y\text{ in }X} \frac{d(fx,fy)}{\alpha\circ d(x,y)}.\] We say $f$ is \emph{$\alpha$--continuous} if $[f]_\alpha<\infty$. More weakly, if each point in $X$ has a neighborhood $U$ such that $\brac*{f\!\restriction_U}_\alpha<\infty$, then   $f$ is said to be $C^\alpha$, or \emph{locally $\alpha$--continuous}. 

Here and throughout this article, we will alternate between the notation $fx$ and $f(x)$ for a function $f$ applied to a point $x$. We will opt
for  the first notation to avoid large numbers of parentheses, and for the second when confusion might arise otherwise.

\begin{lem}\label{l:loc-c-alpha}
If $\alpha$ is a concave modulus,
every compactly supported real-valued locally $\alpha$--continuous map on a metric space is $\alpha$--continuous.
\end{lem}
\bp
This follows from a standard compactness argument, as we illustrate for the convenience of the reader.
Let $f\co X\to \bR$ be the given locally $\alpha$--continuous map.
Pick a finite number of open sets \[U_1,\cdots,U_m,\] where
here \[\overline{\supp{f}} \subset \bigcup_{i=1}^m U_i,\] and $f\restrict_{U_i}$ is $\alpha$--continuous for each $i$.
Set \[U_0 = X\setminus\overline{\supp{f}},\] and let $\delta$ be the Lebesgue number of the open cover \[\{U_0, U_1,\cdots,U_m\}.\]
Then for any distinct $x,y \in X$, we have two cases. If $d(x,y) < \delta$, then $x,y\in U_i$ for some $i$ so that 
\[d(fa,fb) \le \alpha \circ d(x,y) \cdot\max{\brac*{f\!\restriction_{U_i}}_\alpha}.\] If $d(x,y) \ge \delta$, then 
\[d(fa,fb) \le 2\abss{f}_{\infty} \le \alpha \circ d(x,y) \cdot\frac{2\abss{f}_{\infty}}{\alpha(\delta)}.\qedhere\]\ep

Now, suppose that $X$ and $Y$ are smooth, orientable manifolds of dimensions $m$ and $n$, respectively.
For an integer $k\ge0$, a $C^k$--map $g\co X\to Y$ is said to be \emph{$C^{k,\alpha}$}
if each point in $X$ admits a coordinate neighborhood $U$ such that 
the restriction $(D^k g)\restriction_U$ is $\alpha$--continuous.
Here, we regard $U\sse\bR^m$ so that $(D^k g)\restriction_U$ is a map from $U$ to $L((\bR^m)^k,\bR^n)$, the space of $k$--linear maps.

In the sequel, the notation $[\cdot]$ means that the notation in the square bracket forms an optional hypothesis and may be omitted.
For example, if we assert a statement on $C^{k[,\alpha]}$ maps, the statement is meant for $C^k$ and also for $C^{k,\alpha}$ maps.
We will write that a map is $C^{k,\mathrm{bv}}$ if it is $C^k$ and if the $k^{th}$ derivative is of bounded variation.

The \emph{support} of a diffeomorphism $g$ on a manifold $X$ is defined as
\[\supp g = \overline{X\setminus \Fix g} = \overline{\{x\in X\mid gx\ne x\}}.\]
This notion of support is sometimes called the \emph{closed support}, which serves to distinguish it from the
analogously defined \emph{open support} used by some authors. We say $g$ is \emph{compactly supported} if $\supp g$ is compact. We set
\begin{align*}
\Diff_+^{k[,\alpha]}(X)&=\{g\co X\to X\mid g\text{ is orientation preserving }C^{k[,\alpha]}\text{ diffeomorphism}\};\\
\Diff_c^{k[,\alpha]}(X)&=\{g\in \Diff_+^{k[,\alpha]}(X) \mid g\text{ is compactly supported}\};\\
\DD(X;k[,\alpha])&=\{g\in\Diff_c^{k[,\alpha]}(X) \mid \text{there exists an isotopy }\{g_t\}\sse\Diff_c^{k[,\alpha]}(X)\text{ such that }\\
&\qquad\qquad\quad\qquad\qquad g_0=\Id,g_1=g,\text{ and }\cup_{t}\supp g_t\text{ is bounded}\}.
\end{align*} 
All these sets are in fact groups; see Proposition~\ref{prop:diff-group}.

\begin{rem}
 The group $\DD(X;k)$ may not coincide with the group \[\{g\mid g\in \Diff_c^{k}(X), \textrm{ $g$ is isotopic to the identity}\},\]
 even when $M=\bR^n$, with $n\geq 2$. See~\cite{Schweitzer2011ETDS}, for instance.
 \end{rem}

\bd
Let $\alpha$ be a concave modulus.
\be
\item We say $\alpha$ is \emph{sup-tame} if $\lim_{t\to 0+}\sup_{x>0} t\alpha(x) /\alpha(tx)=0$.
\item We say $\alpha$ is \emph{sub-tame} if $\lim_{t\to 0+}\sup_{x>0} \alpha(tx) /\alpha(x)=0$.
\item A \emph{tame pair} is a tuple $(i,\alpha)$ such that one of the following holds:
\begin{itemize}
\item $i=0$ and $\alpha$ is sup-tame;
\item $i=1$ and $\alpha$ is sub-tame;
\item $i\in\bZ\setminus\{0,1\}$.
\end{itemize}
\ee
\ed

For instance, the modulus $\omega_s(x)=x^s$ is both sup-tame and sub-tame for $s\in(0,1)$.
The Lipschitz modulus $\omega_1(x)=x$ is sub-tame, but not sup-tame.

Equipped with this terminology, we posit the following refinements of Mather's results in \cite{Mather1,Mather2}.

\begin{thm}[Mather's Theorem, refined]\label{thm:mather-n}
If $X$ is a smooth $n$--manifold without boundary, and if $(k-n,\alpha)$ is a tame pair,
then the group $\DD(X;k,\alpha)$ is simple.
In particular, if $k\in\bN\setminus\{ n+1\}$ then $\DD(X;k)$ is simple.
\end{thm}

\begin{rem}
The second conclusion of Mather's Theorem holds for $k=0$ and $k=\infty$ as well. The former is classical (Anderson, Chernavski, Kirby, Edwards~\cite{Anderson1958,Rushing1973}) 
and
 the second is due to Thurston (Thurston~\cite{Thurston1974BAMS,Banyaga1997, HT2003,HRT2013,Mann2016NYJM}).
\end{rem}

In this article,  we describe the proof of Mather's Theorem for the case when $n=1$, which reads as follows.
\begin{thm}[Mather's Theorem for one--manifolds, further refined]\label{thm:mather-1d}
The following groups are simple.
\be
\item $\Diff_c^k(\bR)$  and $\Diff_+^{k}(S^1)$, where $k\in\bN\setminus\{2\}$;
\item $\Diff_c^{k,\alpha}(\bR)$ and
$\Diff_+^{k,\alpha}(S^1)$, where $(k-1,\alpha)$ is a tame pair;
\item $\bigcup_{s>r} \Diff_c^s(\bR)$ and $\bigcup_{s>r}\Diff_+^s(S^1)$ for $r\ge1$;
\item $\bigcap_{s<r} \Diff_c^s(\bR)$ and $\bigcap_{s<r}\Diff_+^s(S^1)$ for $r>3$.
\ee
\end{thm}

\begin{rem}
If $k\in\{1,2\}$,
it is not known whether or not $\Diff_c^{k,\alpha}(\bR)$ is simple for all concave moduli $\alpha$.
Let  $\Diff_c^{k+\mathrm{bv}}(M)$ be the group of compactly supported $C^k$ diffeomorphisms $f$ of $M$ such that $\der{f}{k}$ is of bounded variation.
Mather~\cite{Mather3} constructed a nontrivial homomorphism from the group $\Diff_c^{1+\mathrm{bv}}(\bR)$ into $\bR$, thus
 showing that group is \emph{not} simple.
\end{rem}

Let $M\in\{I,S^1,\bR\}$, and let $G\le \Homeo_+(M)$ be a countable group. We define the \emph{critical regularity} of $G$ to be
\[\CR(G)=\sup\{s\mid G\le\Diff_+^s(M)\}.\]
The critical regularity
of a group is, in a precise sense, the smoothest possible incarnation of a group acting on $M$.

The reason for considering countable subgroups of $\Homeo_+(M)$ comes from a result of Deroin--Kleptsyn--Navas~\cite{DKN2007} which
implies that such a group is topologically conjugate to a group of Lipschitz homeomorphisms, so that we generally have the \emph{a priori}
bound $\CR(G)\geq 1$ for any countable group of homeomorphisms.

Critical regularity of groups of homeomorphisms falls into a more general problem of finding
``optimal" regularities $(k,\alpha)$ such that a subgroup
$G\le\Homeo_+(M)$ can be realized as a subgroup of $\Diff_+^{k,\alpha}(M)$, where here $\alpha$ may be a general modulus of continuity
and not just a H\"older modulus. In~\cite{KK2017crit}, the second and third authors developed a general framework for constructing countable
(and in fact finitely generated) groups with specified critical regularities. For details and precise statements of results, we refer the  reader
to Section~\ref{sec:critical}. 

The reason for considering critical regularity in the context of Mather's Theorem is that the methods of~\cite{KK2017crit} can be used
to severely restrict homomorphisms between diffeomorphism groups of $\bR$ for non--integer regularities.

Mann~\cite{Mann2015} proved that there are no injective homomorphisms 
\[
\Diff_c^p(\bR)\to\Diff_c^q(\bR)\]
for integers $q>p\ge 3$. 
From Theorems~\ref{thm:mather-1d} and the critical regularity result~\ref{thm:homomorphism} below,
we can obtain the following generalization of Mann's result.
\begin{cor}\label{cor:cont-hom}
Every homomorphism of the following type is trivial.
\begin{itemize}
\item $\Diff_c^r(\bR)\to\bigcup_{s>r}\Diff_c^s(\bR)$ for each real number $r\in[1,\infty)\setminus\{2\}$;
\item $\bigcap_{s<r}\Diff_c^s(\bR)\to \Diff_c^r(\bR)$ for each real number $r>3$.
\end{itemize}
\end{cor}

\subsection*{Notes and references}
This note is largely based on  Mather's papers on the topic~\cite{Mather1,Mather2}.
Some of the background material follows the original ideas of Higman, Epstein~\cite{Epstein1970} and Ling~\cite{Ling1984}. The reader may refer to a more concisely written account of Thurston's proof for the $C^\infty$ case in~\cite{Banyaga1997}.
Another source of detailed calculations regarding Mather's works for general dimensions can be found in~\cite{HallerThesis}.
Tsuboi gave an alternative proof of Mather's Theorem in~\cite{Tsuboi1985}.

As far as the authors are aware, the generalizations of Mather's Theorem as in parts (2) through (4) of Theorem~\ref{thm:mather-1d}, and as in Corollary~\ref{cor:cont-hom} have not appeared in the literature, except  in~\cite{KK2017crit} where it is stated without proofs. 
Mather considered hypotheses that are more restrictive than the tameness conditions used in this note.
Some of the tools required for our generalizations are necessarily new
and possibly interesting in their own right. We refer the reader to Theorem~\ref{thm:decompose} and to Lemma~\ref{l:tame2}, for instance. 

 Section~\ref{sec:critical} on critical regularity is based on the joint work of the second and the third authors~\cite{KK2017crit}.

\section{Concave moduli}\label{s:mod}

For the remainder of this article, we let $M$ denote a connected one--manifold without boundary; namely, $M\in\{\bR,S^1\}$. 
We will let $k$ denote a positive integer, and let $\alpha\co[0,\infty)\to[0,\infty)$ be a concave modulus.
By convention, we often extend $\alpha$ to $\bR$ by the relation $\alpha(-x)=\alpha(x)$.
We will equip $S^1$ with the usual angular metric $d_{\mathrm{angular}}$.
Note that
  $\DD(\bR;k[,\alpha])=\Diff_c^{k[,\alpha]}(\bR)$
and that
  $\DD(S^1;k[,\alpha])=\Diff_+^{k[,\alpha]}(S^1)$.  

In this section, we establish the following:
\begin{prop}\label{prop:diff-group}
The set $\DD(M;k[,\alpha])$ is a group under composition.
\end{prop}
More importantly, we will prove that each $C^k$--diffeomorphism of $M$ admits an``optimal regularity'' $C^{k,\beta}$ in the following sense. This theorem is particularly useful for us (and seems new to us) since a group written as a union of perfect subgroups is perfect itself.

\begin{thm}\label{thm:decompose}
We have
\begin{align*}
\DD(M;k)&=\bigcup\{\DD(M;k,\beta)\mid
\beta\text{ is a concave modulus}\}\\
&=
\bigcup\{\DD(M;k,\gamma)\mid
\gamma\text{ is a sup-tame modulus}\}.\end{align*}
Moreover, for each family of concave moduli $\FF_0$,
there exists another family $\FF$ of concave moduli
such that 
\[
\bigcap\{\DD(M;k,\mu)\mid
\mu\in\FF_0\}
=
\bigcup\{\DD(M;k,\nu)\mid
\nu\in\FF\}.\]
\end{thm}

\subsection{Group structure}





Let us define a set of integer tuples
\[A_k=\left\{(j_1,\ldots,j_i)\mid 1<i<k, 1\le j_1\le j_2\le\cdots\le j_i, \sum_{t=1}^i j_t=k
\right\}.\]
In particular, we have $A_1 = A_2 = \emptyset$.

\begin{lem}[\cite{Mather1}]\label{l:higher}
If $k \ge 2$, then there exists a set of nonnegative integers  $\{C_\gamma\mid \gamma\in A_k\}$
 such that whenever we are given with $C^k$--maps between open intervals $U,V,W$ as
 \[
\xymatrix{
U\ar[r]^g&V\ar[r]^f& W},\]
we have that
\[
\der{(f\circ g)}{k} =(\der{f}{k}\circ g)(g')^{ k} +(f'\circ g)\der{g}{k} + \sum_{\gamma\in A_k} C_\gamma (\der{f}{i}\circ g)\prod_{t=1}^i \der{g}{j_t}.\]
In particular, if $f$ and $g$ are  $C^{k,\alpha}$, 
then so is $f\circ g$.
\end{lem}
In the above, we used an abbreviation
\[
 \sum_{\gamma\in A_k} \left(\cdot\right)=
\sum_{\stackrel{1<i<k}{\gamma=(j_1,\ldots,j_i)\in A_k}}\left(\cdot\right).\]

\bp The first part is established by a straightforward induction.
Note that the multiplication of two  $C^\alpha$ maps are  $C^\alpha$.
Moreover, the (post- or pre-)composition of a  $C^\alpha$ map with a Lipschitz continuous map is  $C^\alpha$.
The second part thus follows.\ep
\begin{notation}\label{not:pkfg}
We will repeatedly use the last term appearing in the conclusion of Lemma~\ref{l:higher}, which we write more succinctly as 
\[P_k(f,g):=\sum_{\gamma\in A_k} C_\gamma (\der{f}{i}\circ g)\prod_{t=1}^i \der{g}{j_t}.\]
\end{notation}

By setting $g=f^{-1}$ in Lemma~\ref{l:higher}, we obtain the following:
\begin{lem}\label{l:inverse}
For open intervals $U,V\sse M$
and for a $C^1$--diffeomorphism 
$f\co U\to V$, we have that $(f^{-1})'=1/ (f'\circ f^{-1})$.
If $f$ is $C^k$ for some $k\ge2$, then we have
\[
\der{(f^{-1})}{k} =-(\der{f}{k}\circ f^{-1})((f^{-1})')^{k+1} - \sum_{\stackrel{1<i<k}{\gamma=(j_1,\ldots,j_i)\in A_k}} C_\gamma (\der{f}{i}\circ f^{-1})\prod_{t=0}^i \der{(f^{-1})}{j_t},\]
where we set $j_0 =1$. In particular, if $f$ is  $C^{k,\alpha}$ then so is $f^{-1}$.
\end{lem}

Proposition~\ref{prop:diff-group} is immediate from the above two lemmas.

\subsection{Optimal concave moduli}
The concavity of $\alpha$ has several important consequences as follows, which we will repeatedly use. The proof of the next lemma
 is immediate.
\begin{lem}[cf. \cite{KK2017crit}]\label{l:omega}
The following holds for a concave modulus  $\alpha$.
\be
\item
The function $x/\alpha(x)$ is monotone increasing
on $[0,\infty)$.
\item
For all $C>0$ and $x\ge 0$, we have that
\[
\min(C,1)\alpha(x)\le \alpha(Cx)\le\max(C,1)\alpha(x).\]
\ee
\end{lem}

We fix the notation $Tx=x+1$ and consider the group of $1$--periodic homeomorphisms of the real line
\[
\Homeo_\bZ(\bR)=\{f\in \Homeo_+(\bR)\mid f\circ T=T\circ f\}.\]
Observe that the natural projection $p\co \bR\to S^1$ determines a central extension
\[1\to \form{T}\to\Homeo_\bZ(\bR)\to\Homeo_+(S^1)\to1.\]
Here, we identified $S^1=\bR/\bZ$. We will write 
\[\Diff^{k[,\alpha]}_\bZ(\bR)=\Diff_+^{k[,\alpha]}(\bR)\cap\Homeo_\bZ(\bR).\]
We denote by $C_\bZ^{[\alpha]}(\bR)$ the set of real--valued [$C^\alpha$] continuous maps  $f$ on $\bR$ such that $f\circ T =  f$. 

A homeomorphism $f\co\bR\to\bR$ is called \emph{eventually $1$--periodic} if 
$f\circ T(x) = T\circ f(x)$ outside some compact interval.
We write
\[
\Homeo_\mathrm{ep}(\bR)=\{f\in \Homeo_+(\bR)\mid f \text{ is eventually 1-periodic}\},\]
\[\Diff^{k[,\alpha]}_\mathrm{ep}(\bR)=\Diff_+^{k[,\alpha]}(\bR)\cap\Homeo_\mathrm{ep}(\bR).\]
Similarly, we denote by $C_\mathrm{ep}^{[\alpha]}(\bR)$ the set of real--valued  [$C^\alpha$] continuous maps $f$ on $\bR$ such that $f\circ T(x)=f(x)$ for all $x$ outside some compact interval. For $*\in\{c,\bZ,\mathrm{ep}\}$, note that each map in $C_*^{\alpha}(\bR)$ is (globally) $\alpha$--continuous. 

We note that each compactly supported real-valued map has an \emph{optimal} modulus $\beta$:
\begin{prop}\label{prop:opt-conc-mod}
Let $*\in\{c,\bZ,\mathrm{ep}\}$.
Then for each $f\in C_*(\bR)$, there exists a concave modulus $\beta$ such that 
\[f\in C_*^\beta(\bR)=\bigcap \{C_*^\alpha(\bR)\mid \alpha\text{ is a concave modulus and }f\in C_*^\alpha(\bR)\}.\]
Moreover, two such concave moduli $\beta_1,\beta_2$ have the same germs at $0$ up to Lipschitz equivalence.
\end{prop}
The last sentence means that there exist $\delta>0$ and $C>1$ satisfying
\[
1/C \le  \beta_1(x)/\beta_2(x) \le C\]
for all $x\in(0,\delta)$.

\bp
Define $\mu\co [0,\infty)\to [0,\infty)$ by
\[\mu(t):=\sup \left\{|fx-fy| \co  |x-y|\le t\right\}\]
By the uniform continuity of $f$, we see that $\mu$ is monotone increasing, continuous and subadditive.
Therefore, if we define \[\mu_s(t):=(1+t/s)\mu(s),\] we see that $\mu(t)\le\mu_s(t)$ for all $s,t>0$.

We define $\mathcal{F}$ to be the family of continuous, monotone increasing, concave functions \[\nu\co [0,\infty)\to[0,\infty)\]
such that $\mu(x)\le \nu(x)$ for all $x$.
This family is nonempty since $\{\mu_s\co s>0\}\sse\mathcal{F}$.
Define 
\[\beta_0(x):=\inf_{\nu\in\mathcal{F}}\nu(x).\]
Then $\beta_0$ is concave and satisfies $\beta_0(0)=0$.
Moreover, we see that
\[
\mu(x)\le \beta_0(x) \le \inf_{s>0} \mu_s(x)\le\mu_x(x)=2\mu(x).\]
It follows that \[\beta=\beta_0+\Id\co[0,\infty)\to[0,\infty)\]
is a concave modulus satisfying $\mu\le\beta\le 2\mu+\Id$. In particular, we obtain $f\in C^\beta_*(\bR)$.

Now, suppose $f\in C_*^\alpha(\bR)$ for a concave modulus $\alpha$.
We set $R=\diam\supp f$ if $*=c$, and we set $R=1$ otherwise. 
For $0<t\le R$, we have \[ \beta(t) \le 2\mu(t)+t \le\left( 2[f]_\alpha  + R/\alpha(R)\right) \alpha(t).\]
This implies
$C_*^\beta(\bR)\sse C_*^\alpha(\bR)$, whence the desired conclusion.

Finally, for any two concave moduli $\beta_1$, $\beta_2$ satisfying the condition, the last inequality shows that there exists $R > 0$ and $C > 0$ such that for $0 < t \le R$ we have
$$\beta_1 (t) \le C \beta_2(t), \text{          } \beta_2 (t) \le C \beta_1(t).$$
These inequalities imply that they have the same germs at 0 up to Lipschitz equivalence.

\ep
In the proof of Proposition~\ref{prop:opt-conc-mod}, the map $\beta_0$ is called the \emph{least concave majorant (LCM)} of $\mu$.

\bp[Proof of Theorem~\ref{thm:decompose}]
For the case $M=\bR$ we have $\DD(\bR;k[,\alpha])=\Diff_c^{k[,\alpha]}(\bR)$, so that we can substitute $\Diff_c^\cdot(\bR)$ for $\DD(M;\cdot)$. For the case $M=S^1$, we note that the central extension obtained from $p:\bR \to S^1$ gives rise to an exact sequence
\[1\to \form{T} \to \Diff_\bZ^{k[,\alpha]}(\bR) \to \Diff_+^{k[,\alpha]}(S^1) = \DD(S^1;k[,\alpha]) \to 1.\]
Thus, we only need to consider the groups $\Diff_*^{k,[\alpha]}(\bR)$ for $*\in\{c,\bZ\}$.

Let $f\in \Diff_*^{k}(\bR)$. 
By applying Proposition~\ref{prop:opt-conc-mod} to $\der{(f-\Id)}{k}\in C_*(\bR)$,
we obtain a concave modulus $\beta$ such that $f\in\Diff_*^{k,[\beta]}(\bR)$.
Fix $\delta>0$ such that $\beta(\delta)<1$, and define the concave modulus
\[\gamma(x)={\beta}(x)^s\co[0,\infty)\to[0,\infty)\] for $s\in(0,1)$. Then $\gamma(x)$ is sup-tame since $t\beta^s(x)/\beta^s(tx) \le t^{1-s}$ for $x\in(0,\infty)$ and $t\in(0,1)$.
We also have that $0<\beta(x)<\gamma(x)$ for $x\in(0,\delta)$, whence $f\in \Diff_*^{k,\gamma}(\bR)$.

For the second part, $\FF_0$ is given, and we define 
$\FF$ to be the family of concave moduli $\alpha$ such that
\[C_*^\alpha(\bR)\sse \bigcap\{C_*^\mu(\bR)\co \mu\in\FF_0\}.\]
Pick an arbitrary
\[
f\in \bigcap\{
\Diff_*^{k,\mu}(\bR)
 \mid
\mu\in\FF_0\}.\]
By  Proposition~\ref{prop:opt-conc-mod}, there exists a concave modulus $\beta\in\FF$ such that  $ \der{(f-\Id)}{k}\in C_*^\beta(\bR)$.
In particular, $f \in \Diff_*^{k,\beta}(\bR)$.
This implies the conclusion.
\ep

Recall our notation $\omega_s(x)=x^s$, where $s\in(0,1]$.
For the sake of notational compactness, we will write \[\DD(M;k+s)=\DD(M;k,\omega_s).\]
If $r\in\bR$, we denote the largest integer that is strictly less than $r$ as 
\[\form{r}=-\lfloor -r\rfloor-1.\]
For instance $\form{1.1}=1$ and $\form{1}=0$. The second part of Theorem~\ref{thm:decompose} implies the following.

\begin{cor}\label{c:decompose}
Let $r>1$ be real. There exists a family of concave moduli $\FF$ such that
\[\bigcap_{s<r}\DD(M;s)=\bigcup_{\alpha\in\FF}\DD(M;\form{r},\alpha).\]\end{cor}

\subsection{Locally defined concave moduli}
It is often the case that a concave modulus $\mu$ is only locally defined;
that is, $\mu$ is given as a concave homeomorphism $\mu\co [0,\delta)\to[0,\delta')$ for some small $\delta,\delta'>0$.  
For instance, Mather uses moduli which are only defined near zero~\cite{Mather1}. 

In the following lemmas, we observe that such a locally defined concave map  $\mu$ can be replaced by a globally defined concave modulus $\beta\co[0,\infty)\to[0,\infty)$ without affecting the functions which are controlled by it.
Furthermore, the tameness conditions can be certified by the local conditions described below.


\begin{lem}\label{l:tame-criteria}
	Let $\alpha$ be a concave modulus.
	\be
	\item
	The map $\alpha$ is sup-tame if and only if there exists $t_0 > 0 $ satisfying
	\[\sup_{x>0} {t_0\alpha(x)}/{\alpha(t_0 x)} < 1.\]
	\item
	The map $\alpha$ is sub-tame if and only if there exists $t_0 > 0 $ satisfying
	\[\sup_{x>0} {\alpha(t_0 x)}/{\alpha(x)} < 1.\]
	\ee
\end{lem}
\bp
Only if part is clear from the definition of tameness.

\item
For (1), suppose that $t_0 >0$ satisfies the condition of the lemma. 
Since $x/\alpha(x)$ is monotone increasing, so is the map
 \[F(t) := \sup_{x>0} {t\alpha(x)}/{\alpha(tx)}=t\sup_{x>0} {\alpha(x)}/{\alpha(tx)},\] 
which is defined for each $0<t\le t_0$.
Moreover, for each $0<t \le \min(t_0 , \sqrt{t_0})$ we have
\[F(t^2)=\sup_{x>0}t^2\alpha(x)/\alpha(t^2x) =\sup_{x>0} t\alpha(x)/\alpha(tx) \cdot t\alpha(tx)/\alpha(t^2x)\le F(t)^2.\]

Since the function $F$ is monotone increasing and nonnegative, the limit
$\lim_{t\to 0+} F(t)$ exists. For all positive integer $n$ we have
\[
0\le
\lim_{t\to 0+} F(t)
\le 
F(t_0^{2^n}) \le F(t_0)^{2^n}.
\]
Hence, the desired limit is 0.
The proof of part (2) is nearly identical after replacing $F(t)$ by 
\[G(t) = \sup_{x>0} {\alpha(tx)}/\alpha(x).\qedhere\]
\ep

The following lemma is useful for the flexibility which it furnishes.
First, one may require a  concave moduli only to be defined near zero without losing generality.
Second, the tameness is also determined by the local behavior near zero, which is often simpler to verify.
Recall two maps $f,g$ defined near $0$ are said to have the same \emph{germs} at $0$
if $f(x)=g(x)$ for all $x$ that are sufficiently near $0$.

\begin{lem}\label{l:tame2}
	For each concave homeomorphism $\alpha\co [0,\delta)\to[0,\alpha(\delta))$,
	there exists a concave modulus $\beta\co[0,\infty)\to[0,\infty)$ such that $\alpha$ and $\beta$ have the same germs at $0$.
	Moreover, the following hold.
	\be
	\item
	The map $\beta$ can be chosen to be a sup-tame modulus 
	if
	\[\lim_{t\to0+}\sup_{x<\delta} {t\alpha(x)}/{\alpha(tx)}=0.\]
	\item
	The map $\beta$ can be chosen to be a sub-tame modulus
	if \[\lim_{t\to0+}\sup_{x<\delta} {\alpha(tx)}/{\alpha(x)}=0.\]
	\ee\end{lem}
\bp
For the first claim, we will extend $\alpha \restrict_{[0,\delta/2]}$ to a concave modulus $\beta$ of the form \[\beta(x) = a\sqrt{x}+b \quad\quad (a,b \in \bR, \quad a>0)\] for $x \ge \delta/2$.
Here, the constants $a$ and $b$ are determined by the conditions
$\beta(\delta/2)=\alpha(\delta/2)$
and
\[
\beta'(\delta/2)=\frac{\alpha(\delta)-\alpha(\delta/2)}{\delta/2}>0.\]
It is clear that $\beta\co[0,\infty)\to[0,\infty)$ is a homeomorphism. Also, $\beta$ is concave since the restrictions $\beta \restrict_{[0,\delta/2]}$ and $\beta \restrict_{[\delta/2, \infty)}$ are both concave, and since the line joining $(x,\beta(x))$ to $(y,\beta(y))$ for $x<\delta/2<y$ is below the line passing through the point $(\delta/2 , \beta(\delta/2))$ with the slope 
\[\frac{\alpha(\delta)-\alpha(\delta/2)}{\delta/2}.\]


It remains to show parts (1) and (2). For part (1), let us pick $0<t_0<1$ such that 
\[\sup_{x<\delta} {t_0\alpha(x)}/{\alpha(t_0 x)}<1.\]
Since $\beta$ is nonlinear and concave on $[\delta,\infty)$, we see that 
\[{t_0\beta(x)}/{\beta(t_0 x)}< 1\] for $x>\delta$ as well.
From \[\lim_{x\to \infty} {t_0\beta(x)}/{\beta(t_0 x)} = \sqrt{t_0}<1,\]
we see that 
 $\sup_{x>0} {t_0\beta(x)}/{\beta(t_0 x)} <1$.
Lemma~\ref{l:tame-criteria} (1) then implies that $\beta$ is actually sup-tame.

The same construction also works for the second claim, combined with Lemma~\ref{l:tame-criteria} (2).
\ep


We let $<_\bC$ denote the lexicographical ordering on $\bC$, viewed as pairs $(a,b)$ of real numbers with the identification $(a,b)=a+bi$.  
The intervals $(a,b)_\bC$ and $(a,b]_\bC$ are defined by this ordering.
We will repeatedly use the following modulus $\omega_z$ in this paper.

\begin{lem}\label{l:omegaz}
For each complex number  $z=\sigma+\tau i\in(0,1]_\bC$, there exists a positive number $\delta$ such that 
the map \[\omega_z(x):=\exp\left(-\sigma\log(1/x)-\tau\frac{\log(1/x)}{\log\log(1/x)}\right)\]
defined on $(0,\delta]$ extends to a concave modulus on $[0,\infty)$. 
Furthermore, such a concave modulus (still denoted as $\omega_z$) is sup-tame when $\repart z<1$, and sub-tame when $\repart z>0$.
\end{lem}
We emphasize again that such an extension may not be unique, but they all have the same germs. In particular, the notation $\omega_z(x)$ is not ambiguous for all small $x>0$. Note that if $\sigma\in(0,1)$, then we have the H\"older modulus 
\[
\omega_\sigma(x)=x^\sigma.\]

\bp[Proof of Lemma~\ref{l:omegaz}]
We may assume $z\ne1$. 
Let $x$ be sufficiently small. 
Substituting \[L = \log \log(1/x)=\log(-\log x)\quad\text{ and }\quad
f = \log x (\sigma + \tau/L),\]  we have $\omega=\exp\circ f$. It follows that
\[
\omega''=  \left(f''+(f')^2\right)\omega .\]
Using $L' = 1/(x\log x)$, we see that
\begin{align*}
f' &= \frac1x\left( \sigma+\frac\tau L-\frac \tau {L^2}\right),\\
f'' &= -\frac\sigma{x^2}
-\frac{\tau}{x^2L}\left(
1-\frac1L +\frac1{L\log x}-\frac2{L^2\log x}\right).
\end{align*}

We obtain the following.
\begin{align*}
\omega''/\omega 
&= \frac{\sigma^2-\sigma}{x^2} +  O \left(\frac1{x^2L}\right)&\text{if }0<\sigma<1,\\
\omega''/\omega &= -\frac{\tau}{x^2L} +  O \left(\frac1{x^2L^2}\right)&\text{if }\sigma=0\text{ and }\tau>0,\\
\omega''/\omega &= \frac{\tau}{x^2L} +  O \left(\frac1{x^2L^2}\right)&\text{if }\sigma=1\text{ and }\tau<0.
\end{align*}

Hence, there exists $\delta > 0$ such that $\omega''(x) < 0$ for $0 < x \le \delta$.
By Lemma~\ref{l:tame2}, the map $\omega_z \restrict_{[0,\delta]}$ extends to a concave modulus.

Let us now check the tameness conditions.

\begin{claim*}
Let  $z \in (0,1)_\bC$. 
Then $\omega_z$ is sup-tame if and only if $\omega_{1-z}$ is sub-tame.
\end{claim*}
This claim is a consequence of the following equation.
 \[t\omega_z(x)/\omega_z(tx) = \omega_{1-z}(tx)/\omega_{1-z}(x).\]

\begin{claim*}
Let $0<_\bC w <_\bC z <_\bC 1$. If $\omega_z$ is sup-tame, then so is $\omega_w$.\end{claim*}

The map $\omega_z / \omega_w= \omega_{z-w}$ is a concave modulus (and hence, increasing) near $0$.
So, for all $x$ near $0$ and for all $t\in(0,1)$ we have
\[\frac{t\omega_z(x)}{\omega_z(tx)} > \frac{t\omega_w(x)}{\omega_w(tx)} > 0.\]
This implies the claim.

Now, for the case $\sigma = \repart z < 1$, we have $z <_\bC \frac{1+\sigma}{2}$. Since $\omega_{(1+\sigma)/2}$ is sup-tame, so is $\omega_z$ by the second claim.
The desired conclusion for the case $\repart z>0$ follows from the first claim.
\ep
\begin{rem}
If $z=1+\tau i$ for some $\tau<0$, then 
for $g(x) = x/\log(x)$ we have that
\[ \frac{t\omega_z(x)}{\omega_z(tx)} = \exp \left[ -\tau \left( g(\log(1/x)) - g(\log(1/x) + \log(1/t)) \right) \right]. \]
Since $\lim_{x \to \infty} g'(x) = 0$, the value of $t\omega_z(x)/\omega_z(tx)$ converges to 1 as $x \to 0$ for all fixed $0<t<1$. This proves that $\omega_z$ is not sup-tame in this case. 
\end{rem}

\section{Simplicity of commutator groups}\label{s:higman}
Recall our convention that $k$ denotes a positive integer and $\alpha$ denotes a concave modulus.
In this section, we establish the following.
\begin{thm}\label{t:r-comm}
The commutator group of $\DD(M;k[,\alpha])$ is simple.
\end{thm}
We will use Higman's Theorem for the case $M=\bR$,
and Epstein--Ling Theorem for $M=S^1$, as described below.

\subsection{Higman's Theorem}
The case $M=\bR$ of Theorem~\ref{t:r-comm} follows from a classical technique due to Higman.
Consider a group $H$ acting on a set $X$. For $g\in H$, we write \[\supt g=\{x\in X\mid g(x)\neq x\}.\]
Observe that if $g$ is a homeomorphism on a space $X$, then we have
\[
\supp g=\overline{\supt  g}.\]
For $g,h\in H$, we use the notation $g^h = h^{-1}gh$. 
We denote by $\fform{t}_H$ the normal closure of $t\in H$.

\begin{thm}[Higman]\label{thm:higman}
Let $H$ be a group acting faithfully on a set $X$. Suppose that for all triples $r,s,t\in H\setminus \{1\}$ there is an element $u\in H$ such that \[u(\supt  r\cup\supt  s)\cap tu(\supt  r\cup\supt  s)=\varnothing.\] 
Then the following hold.
\be
\item
The commutator subgroup $H'=[H,H]$ is simple and nonabelian.
\item Every proper quotient of $H$ is abelian.
\item $H$ has trivial center.
\ee
\end{thm} 
We will include a proof of Higman's Theorem, as the proof of this version is not very easy to find in the literature;
the reader may also consult \cite{BurilloBook}.

\bp[Proof of Higman's Theorem]
We may certainly assume $H$ is nontrivial.
\begin{claim*}
The following hold:
\be
\item $H'$ is nonabelian.
\item For all $r,s,t\in H\setminus1$, there exists some $u\in H$ such that
\[\brac{r^{t^u},s}=1.\]
\item\label{p:normal-hprime}
Every nontrivial normal subgroup of $H$ contains $H'=[H,H]$.
\item For all $r,s,t\in H\setminus1$, there exists some $w\in H'$ such that 
\[w(\supt  r\cup\supt  s)\cap tw(\supt  r\cup\supt  s)=\varnothing.\]
\ee\end{claim*}

By setting \[r=s=t\in H\setminus\{1\}\] and applying the hypothesis, we see that $H$ is nonabelian.
Then we let \[r=s=t\in H'\setminus\{1\}\] and see that \[r \ne r^u \cdot r \cdot (r^{u})^{-1}\] for some $u\in H$.
Since $r^u\in H'$, we have part (1).

Part (2) follows from the existence of $u\in H$ such that \[\supt  r\cap t^u\supt  s=\varnothing.\] 

For part (3), let $1\ne t\in N\unlhd H$ and pick arbitrary $r,s\in H\setminus1$.  
By part (2), we can find $u\in H$ such that 
\[[r,s]=\left[r\cdot (t^{u})^{-1}r^{-1}t^u,s\right]\in\fform{t}_H\le N.\]
This means $H'\le N$.

In order to prove part (4), we let $u$ be given by applying the hypothesis to the triple $(r,s,t)$.
We may assume $u\not\in H'$; in particular, $u\ne1$.
Applying the hypothesis again to the triple $(t,u,t)$, we obtain $v\in H$ such that 
\[\left[t^{t^v},u\right]=1.\]
Put $w:=\left[t^v,u^{-1}\right]\in H'$ so that 
\[
[t,wu^{-1}]=\left[t,t^v u^{-1}(t^{v})^{-1}\right]=1.\]
Then, we have $t^w = t^u$, so that $w$ is the desired element of $H'$.

We can now finish the proof of the theorem.
Since the center $Z(H)$ is an abelian normal subgroup of $H$, parts (1) and (3) above imply that $Z(H)=1$.

It only remains to show that $H'$ is simple. 
Suppose we have
\[1\ne t\in N\unlhd H'.\]
Let us consider arbitrary $r,s\in H'\setminus1$.
From part (4) above, we have some
 $w\in H'$ such that
\[\brac*{r^{t^w},s}=1.\]
Then we have that
\[
[r,s]=\brac*{r\cdot (t^{w})^{-1}rt^w,s}\in \fform{t}_{H'}\le N.\]
This implies that $H''\le N$. 
Since $1\ne H''\unlhd H$, we see from part (\ref{p:normal-hprime}) that $H'= H''$.  We conclude that $H'\le N$, as desired.
\ep

The following topological variation of Higman's Theorem is particularly useful for our setting.

\bd[\cite{KKL2017}]\label{defn:cpt-trans}
We say $H\le\Homeo(X)$ acts \emph{CO-transitively} (that is,
\emph{compact--open-transitively}) on a topological space $X$ if for each proper compact subset $A\sse X$ and for each nonempty open subset $B\sse X$, there is $u\in H$ such that $u(A)\sse B$.
\ed

\begin{rem}
In the case when $X=S^1$, the CO--transitivity of a group action is equivalent to being ``minimal and strongly proximal''~\cite{Furstenberg1963AM}.
This equivalence is noted in~\cite[Definition 3.4]{BFS2006}, where it is also proved that every non-elementary (non-measure preserving) topological circle action of a group is semi-conjugate to a CO--transitive action.  
\end{rem}


\begin{lem}[\cite{KKL2017}]\label{l:co-trans}
Let $H$ is a group of compactly supported homeomorphisms on a non-compact Hausdorff space $X$ such that $H$ acts
CO--transitively on $X$.
Then $H'$ is simple and nonabelian.
Moreover, every proper quotient of $H$ is abelian, and the center of $H$ is trivial.
\end{lem}

\bp
It is obvious that $H$ is nontrivial.
In order to apply Higman's Theorem, we fix $r,s,t\in H\setminus\{1\}$.
Since $X$ is Hausdorff, we can choose a nonempty open set $B\sse X$ such that $B\cap tB=\varnothing$.
Choose a proper compact subset $A$ containing $\supt  r\cup\supt  s$.
By CO-transitivity, there is $u\in H$ such that $uA\sse B$.
From $uA\cap tuA=\varnothing$, we see
\[
u(\supt  r\cup \supt  s)\cap tu(\supt  r\cup \supt  s)=\varnothing,\]
as desired in the hypotheses of Higman's Theorem.
\ep

Since $\Diff_c^\infty(\bR)$ acts CO--transitively on $\bR$, we obtain the following.
\begin{lem}
If $G$ is a group of compactly supported homeomorphisms on $\bR$ such that 
$\Diff_c^\infty(\bR)\le G$,
then  $G'$ is simple and every proper quotient of $G$ is abelian.
\end{lem}
The case $M=\bR$ of Theorem~\ref{t:r-comm} is now immediate.


\subsection{Fragmenting diffeomorphisms}
Let $X$ be a space and let $G\le\Homeo(X)$. 
For a (usually, open) subset $U\sse X$, we write
\[
G[U]:=\{ g\in G\mid \supp g\sse U\}.\]
We say the action $G\curvearrowright X$ is \emph{fragmented} if for every open cover $\UU$ of $X$ we have 
\[G=\form{ G[U]\mid U\in\UU}.\]
In other words, every $g\in G$ can be written as 
\[g = g_1g_2\cdots g_\ell\]
for some some $\ell\ge0$, $U_1,\ldots,U_\ell\in \UU$ and $g_i\in G[U_i]$.

\begin{lem}\label{l:fragmented}
The group $\DD(M;k[,\alpha])$ is fragmented.
\end{lem}

We will prove the case $M=S^1$,
which is what we require in this article. We remark that the case $M=\bR$ is slightly simpler, though
 almost identical.

\begin{lem}\label{l:top-gp-gen}
Let $G$ be a group equipped with a topology, and suppose that $G$ is connected in this topology.
If the right multiplication $R_g\co G\to G$ is continuous for each $g\in G$, then every neighborhood of the identity generates $G$.\end{lem}
We note that $G$ is not necessarily required to be a topological group.
\bp[Proof of Lemma~\ref{l:top-gp-gen}]
The map $R_g$ is a homeomorphism for each $g\in G$, since $R_g\circ R_{g^{-1}}=\Id$.
Let $W$ be an open neighborhood of the identity, and let $H=\form{W}$.
If $h\in H$, then $Wh=R_h(W)$ is an open neighborhood of $h$ contained in $H$. 
If $g\in G\setminus H$, then $Wg=R_g(W)$ does not intersect $H$; for otherwise, we would have
\[g\in \form{W^{-1},H}=H.\]
In particular, $G\setminus H$ is open. Since $G$ is connected and $H$ is clopen, we see that $G=H$.\ep

For the purpose of this subsection, we will employ the $C^1$ topology on diffeomorphism groups.
Let $\|\cdot\|_\infty$ denote the $L^\infty$--norm.
We define the \emph{$C^1$--metric} on $f,g\in \DD(S^1;1)$ as
 \[d_1(f,g):=\|f-g\|_\infty+\|f'-g'\|_\infty.\]
We recall again that \[\|f-g\|_\infty:=\sup_{x\in S^1} d_{\mathrm{angular}}(fx,gx).\]


\begin{lem}\label{l:top-group} 
The following hold.
\be
\item
The group $(\DD(S^1;1),d_1)$ is a connected topological group.
\item
The subgroup $\DD(S^1;k[,\alpha])\le \DD(S^1;1)$ is connected.
\ee
\end{lem}
\bp
Let $a,b\in\DD(S^1;1)$, and let  $\delta>0$.
Suppose $f,g\in\DD(S^1;1)$ such that \[d_1(a,f),d_1(b,g)<\delta.\]
Then it is routine to see that if $\delta$ is sufficiently small, then $d_1(a^{-1},f^{-1})$ and $d_1(a\circ b,f\circ g)$ are arbitrarily small.
Thus we see that $(\DD(S^1;1),d_1)$ is a topological group.

It only remains to prove the connectivity of $\DD(S^1;1)$ and $\DD(S^1;k[,\alpha])$ with the  $d_1$--topology. 
Suppose $g\in \DD(S^1;k[,\alpha])$. We claim that there is an isotopy
\[H=H(x,t)=H_t(x)\co S^1\times I\to S^1\]
such that the following hold:
\be[(i)]
\item $H_0(x)=x$ and $H_1(x)=g(x)$;
\item For all $\epsilon>0$, there exists a $\delta>0$ such that  whenever $|s-t|<\delta$, 
we have $d_1(H_s,H_t)<\epsilon$;
\item  $H_s\in\DD(S^1;k[,\alpha])$ for each $s\in I$.
\ee

To see this claim, we first apply a suitable rotational isotopy 
$H_t(x)=x+t$
so that we may assume $g(0)=0$. We may thus cut open $S^1$ at the fixed point and view $g$ as a map $g\co I\to I$.
Define \[H_t(x)=tg(x) + (1-t)x.\]
Note that 
\begin{align*}
H_t'(x)&=tg'(x)+(1-t)=1+t(g'(x)-1)>0,\\
|H_s(x)-H_t(x)|&\le |s-t|\cdot \|g\|_\infty + |s-t|\cdot x,\\
|H_s'(x)-H_t'(x)|&\le |s-t|\cdot \|g'\|_\infty.\end{align*}
The first of these inequalities proves that for each $t$, the map $H_t(x)$ is indeed a homeomorphism of $I$.
It follows that $H_t$ is a $C^{k[,\alpha]}$ isotopy,
and that $\|H_s-H_t\|_1<\epsilon$ for all small $|s-t|$.
\ep

For a group $G$ and a subset $S\sse G$, the \emph{subgroup generated by $S$} means
\[\form{S}:=\{ s_1^{e_1}\cdots s_\ell^{e_\ell}\mid \ell\ge0, \ s_i\in S,\ e_i=\pm1\}.\]
As before, we denote by $R_g$ the right-multiplication map $R_g(f)=fg$ for $g\in G$.

\begin{lem}\label{l:frag}
The action of $\DD(S^1;k[,\alpha])$ on $S^1$ is fragmented.\end{lem}
\bp

Let $\VV=\{V_1,\ldots,V_m\}$ be a finite cover of $S^1$, and
let \[\{\phi_i\co S^1\to[0,1]\co 1\le i\le m\}\] be a smooth partition of unity subordinate to $\VV$.
We put \[K=2+\sum_i\|\phi_i'\|_\infty,\]
and choose $0<\epsilon<1/(2K)$.
Consider the following $C^1$--open neighborhood of the identity map:
\[
W=\{g\in \DD(S^1;1)\co
\|g-\Id\|_\infty<\epsilon\text{ and }
\|g'-1\|_\infty<\epsilon\}.\]
By Lemma~\ref{l:top-group}, the group $\DD(S^1;k[,\alpha])$ is a connected subgroup of the topological group $\DD(S^1;1)$.
As $W_0:=W\cap\DD(S^1;k[,\alpha])$ generates the topological group $\DD(S^1;k[,\alpha])$ by Lemma~\ref{l:top-gp-gen},
it suffices to show that each element of $W_0$ is fragmented for $\VV$.

Let $g\in W_0$. We define $g_0=\Id$ and set
\[
g_j(x) := x + \sum_{i=1}^j \phi_i(x)\left(g(x)-x\right)\]
for $1\le j\le m$. Note that $\|g_j-\Id\|_\infty<1/2$.
By the choice of $\epsilon$, we see that 
\[g_j' =1 + (g-x)\sum\phi_i' + (g'-1)\sum\phi_i
\ge 1 - \epsilon\sum_i \|\phi_i'\|_\infty- \epsilon\cdot 1>0,\]
so that we may conclude $g_j\in \DD(S^1;k[,\alpha])$. 
For $1\leq j\leq m$,  if $x$ does not belong to the interior of $V_j$, then $g_{j-1}(x)=g_j(x)$.
In other words, we have $\supp (g_{j-1}^{-1}g_j)\sse V_j$.
The following product provides the desired fragmentation of $g$:
\[
g = (g_0^{-1}g_1)(g_1^{-1}g_2)\cdots (g_{m-1}^{-1}g_m).\qedhere\]
\ep

\subsection{The Epstein--Ling Theorem}

Recall that a Hausdorff space is \emph{paracompact} 
if and only if every open cover $\VV$ admits an open refinement $\UU$
such that  for all sets $A,B\in\UU$, we have that either $A\cap B=\varnothing$, or the union
$A\cup B$ is contained in a third set $C\in\VV$.
Such a refinement $\UU$ is called an \emph{open star refinement} of $\VV$. 
For us, a smooth manifold is always paracompact and admits a smooth partition of unity; see~\cite{Willard2004book} for more details. 

 A \emph{basis} for a topological space is a collection of open sets such that every open set is the union of elements of the basis.
Let  $X$ be a topological space, and let $G\le \Homeo(X)$.
The action $G\curvearrowright X$ is \emph{transitive--inclusive} (with respect to some basis $\BB$)
if for each $U,V\in\BB$, there exists an element $g\in G$ such that $gU\sse V$.

The following result is a variation by Ling of a result originally due to Epstein, which in turn is very much in the same spirit as Higman's Theorem.
\begin{thm}[{Epstein--Ling~\cite{Epstein1970,Ling1984}}]\label{thm:epstein-ling} If a group $G$ admits a faithful, fragmented, transitive--inclusive action on some paracompact space, then $G'$ is
simple.
\end{thm}
\begin{rem}
It will follow from the proof that
every proper quotient of $G$ is abelian and $G$ has trivial center.\end{rem}

The remainder of this subsection is devoted to proving Theorem~\ref{thm:epstein-ling}.
Our proof is based on the ideas of Epstein, Ling and Michalik--Rybicki~\cite{MR2011TA}.

Let $X$ be a paracompact topological space, let $\BB$ be a basis for the topology of $X$, let $G$ be a nonabelian group, and let $G\le\Homeo(X)$ be a faithful fragmented action which is transitive--inclusive with respect to $\BB$. 
Our goal is to show that $G'$ is simple.

For a group $H$, the symbol $H^{(n)}$ denotes the $n$--th term in the solvable series of $H$. 

\begin{lem}\label{l:ling-nonabel}
Let $U$ be an open neighborhood of $x\in X$.  Then the following hold.
\be
\item
There exists $g\in G[U]$ such that $gx\ne x$;
\item
For all integer $n\ge0$,
there exists an element $g\in G[U]^{(n)}$ such that $gx\ne x$. 
In particular, $G[U]$ is not solvable. 
\ee
\end{lem}
\bp
(1)  From the transitive--inclusivity, we see that $G$ does not have a global fixed point. In particular, there exists $a\in G$ such that $ax\ne x$. By fragmenting $a$ with respect to the cover $\{U,X\setminus\{x\}\}$, we find some $g\in G[U]$ such that $gx\ne x$.

(2) 
Suppose we have shown the claim for the case $n\ge0$.
We have some $g\in G[U]^{(n)}$ such that $gx\ne x$. 
Choose an open neighborhood $V\sse U$ of $x$ such that $gV\cap V=\varnothing$. 
Applying the inductive hypothesis to the pair $(x,V)$, we can find some $h\in G[V]^{(n)}$ such that $hx \ne x$. Then we have that 
\[ghx\ne gx=hgx.\]
Then $[g^{-1},h^{-1}]\in G[U]^{(n+1)}$ is the desired element for $n+1$.
\ep
For each $x\in X$ there exists some open neighborhood $V_x\in\BB$ of $x$ and some $a\in G$
such that $V_x\cap a(V_x)=\varnothing$.
We define an open cover \[\VV :=\{ V_x\mid x\in X\}\]
and let $\UU$ be an open star refinement of $\VV$.
\begin{lem}\label{l:ling-prop-abel}
The following hold.
\be
\item
Every proper quotient of $G$ is abelian.
\item
$G'=G''$.
\item
$G$ has trivial center.
\ee
\end{lem}
\bp
(1)
Suppose we have $1\ne t\in N\unlhd G$.
Pick an element $V(t)\in \BB$ such that $V(t)\cap tV(t)=\varnothing$.
We consider an arbitrary element $V\in \VV$
and let $f,g\in G[V]$. There exists some $u\in G$ such that $uV\sse V(t)$, so that  \[u\supp f\cap tu\supp g\sse uV\cap tuV=\varnothing.\]
Put $v:=t^u (g^{-1}) (t^u)^{-1}$. Then we have $[f,v]=1$ and that 
\[[f,g]=[f,gv]=\left[f,\left[g,t^u\right]\right]\in \fform{t}_G.\]
It follows that $G[V]'\le  \fform{t}_G \le N$.

Now let $a,b\in G$ be arbitrary.
Applying fragmentation, we can write 
\[
a = \prod a_i,\quad b=\prod b_j\]
such that $\supp a_i\sse U'_i, \supp b_j\sse U''_j$ for some suitable elements $U_i', U_j''\in \UU$.
We have the commutator formula
\[
[fg,h]=[f,h]\cdot  [g,h]^u,\quad [f,gh]=[f,g]\cdot[f,h]^v\]
for some words $u,v$ in $\{f,g,h\}$.
So, we can write 
\[
[a,b]=\prod_{i,j} [a_i,b_j]^{u_{ij}}\]
for some suitable words $u_{ij}$ in the group elements $\{a_1,\ldots,b_1,\ldots\}$.
Moreover, the $(i,j)$--term in the product is nontrivial only if 
\[\varnothing\ne \supp a_i\cap \supp b_j\sse U_i'\cap U_j''.\]
For such $(i,j)$ we can find $V_{ij}\in\VV$ such that 
$a_i,b_j\in G[V_{ij}]$.
From the previous paragraph, we see that $[a,b]\in N$.
It follows that $G'\le N$, as desired.

(2) This is immediate from the fact that $1\ne G''\unlhd G$ and from part (1).

(3) Note that the center $Z(G)$ is an abelian normal subgroup of $G$. 
Since $G'=G''$ is nonabelian, we see $G'\not\le Z(G)$.
By part (1), we conclude that $Z(G)=1$.
\ep

\begin{lem}\label{l:comm-replace} 
For each  $u\in G$ and for each $U\in \UU$ there exists $w\in G'$ such that $uU=wU$. \end{lem}
\bp
Since $G$ is fragmented for $\UU$, we can find 
$U_1,\ldots,U_\ell\in\UU$ such that \[ u = u_\ell \cdots u_1\] for some $u_i\in G[U_i]$. 

If $U_i\cap U=\varnothing$, we set $a_i=1$. 
If $U_i\cap U\ne\varnothing$, then there exists $V\in \VV$ such that $U_i\cup U\sse V$; in this case, we can find
an element $a_i\in G$ such that $ U \cap a_i U_i=\varnothing$.
For $i>0$, we define
\[
g_i = (u_{i-1}\cdots u_1) a_i\in G.\]
Then we see that 
\[ g_i U_i\cap \left(u_{i-1}\cdots u_1U\right)=\varnothing.\]
Since
$\supp(g_i u_i g_i^{-1})$ is disjoint from $u_{i-1}\cdots u_1U$,
we see that
\begin{align*}
\prod_{i=\ell}^1[u_i,g_i]U
&=\prod_{i=\ell}^2[u_i,g_i]\cdot u_1\cdot g_1 u_1^{-1} g_1^{-1}U
=\prod_{i=\ell}^2[u_i,g_i]\cdot u_1U\\
&=\prod_{i=\ell}^3[u_i,g_i]\cdot u_2\cdot g_2 u_2^{-1} g_2^{-1}\cdot u_1U
=\prod_{i=\ell}^3[u_i,g_i]\cdot u_2u_1U=\cdots=uU.\qedhere
\end{align*}
\ep
\bp[Proof of Theorem~\ref{thm:epstein-ling} for $M=S^1$]
Let us fix an arbitrary $1\ne s\in G'$.
It suffices for us to show that \[G'=\fform{s}_{G'}.\]

By the preceding lemmas, we can define
\[
1\ne H:=\left\langle \left\{ G[uU]''\mid u\in G, U\in \UU \right\}\right\rangle
=\left\langle \left\{ G[wU]''\mid w\in G', U\in \UU \right\}\right\rangle\le G''=G'.\]
For each $u,g\in G$ and $U\in\UU$
we have
\[gG[uU]''g^{-1}=G[guU]''\le H.\]
This implies that $H\unlhd G$. 
By Lemma~\ref{l:ling-prop-abel}, we have $H=G'$.

We have some nonempty set $V(s)\in \BB$ such that $V(s)\cap sV(s)=\varnothing$. 
Pick an arbitrary set $U\in\UU$ and elements $f,g\in G[U]'$.
By the transitive--inclusivity, we have some $u\in G$ such that $uU\sse V(s)$.
From Lemma~\ref{l:comm-replace}, we can find $w\in G'$ with $wU=uU$.
Since \[w\supp f\cap sw\supp g\sse wU\cap swU=\varnothing,\]
we have that
\[
[f,g]=\left[f,\left[g,s^w\right]\right]\in  \fform{s}_{G'}.\]
This implies that $G[U]''\le \fform{s}_{G'}$.

For all $w\in G'$, we have
\[G[wU]'' = wG[U]'' w^{-1}\le  \fform{s}_{G'}.\]
From the definition of $H$, we see that $H=G'\le \fform{s}_{G'}$.
This implies $G' = \fform{s}_{G'}$, as desired.
\ep

\bp[Proof of Theorem~\ref{t:r-comm} for $M=S^1$]
We have seen that \[\DD(\bR;k[,\alpha])\] is fragmented.
Since $\DD(\bR;\infty)$ is transitive--inclusive, we deduce the result from the Epstein--Ling Theorem.
\ep

\subsection{The piecewise--linear case}

Let $\PL(I)$ denote the group of orientation--preserving, piecewise linear homeomorphisms of $I$ with finitely many break points. 
\begin{cor}
The double commutator subgroup $\PL(I)''$ is simple.\end{cor}
\bp
Consider the ``germ homomorphism''
\[
\Phi\co \PL(I)\to \bR^2\]
defined by $\Phi(f) = (\log f'(0),\log f'(1))$. The map $\Phi$ is a group homomorphism onto the abelian group  $\bR^2$. In particular,  $ \PL(I)'\le\ker\Phi$
and $\PL(I)'\le \Homeo_c(0,1)$.

It suffices to show that $\PL(I)'$ is CO--transitive on $(0,1)$, by Lemma~\ref{l:co-trans}.
Let $A\sse (0,1)$ be compact,
and let $B\sse (0,1)$ be open. Then we can find
an element $f\in \PL_c(0,1)$ such that $fA\sse B$, where $\PL_c$ denotes the subgroup
consisting of homeomorphisms supported compactly in $(0,1)$.
Let $C\sse(0,1)$ be a compact interval containing  $A\cup \supp f$. There exists an element $h\in \PL_c(0,1)$ such that
\[\supp hf^{-1}h^{-1}\cap\supp f \sse hC\cap C=\varnothing.\]
We see that 
\[[f,h](A) = f\cdot hf^{-1}h^{-1}(A)=f(A)\sse B.\]
This implies that $\PL(I)'$ acts CO--transitively on $(0,1)$.
\ep

\begin{cor}
The group $[\PL_c(I),\PL_c(I)]$ is simple.\end{cor}
\bp
Note $\PL(I)'\le  \PL_c(I)=\ker\Phi$, which implies that $\PL_c(I)$ is also CO--transitive.
\ep


\section{Eventually periodic diffeomorphisms}
In this section, we study the $C^{k[,\alpha]}$--topology of diffeomorphism groups. 
There are three main parts to this discussion. First, we investigate how this topology is related to the group structure.
Second, we estimate the behavior of the $C^{k,\alpha}$--norm under composition of functions.
Third, we observe that a sufficiently small closed $C^{k,\alpha}$--neighborhood of the identity enjoys the fixed point property.

\subsection{The $C^{k}$ and $C^{k,\alpha}$ topologies}\label{ss:notations}
For each $\beta\in\bR$, we let $T_\beta$ denote the translation $T(\beta)(x)=T_\beta(x)=x+\beta$,
and we fix the notation $T:=T_1$.  For a compact interval $J\sse\bR$, we put
\begin{align*}
\Homeo_J(\bR)&=\{f\in\Homeo_+(\bR)\mid \supp f\sse J\};\\
\Homeo_\bZ(\bR)&=\{f\in \Homeo_+(\bR)\mid f\circ T = T\circ f\}.
\end{align*}

\begin{notation}\label{not:If}
For a homeomorphism $f$ of $\bR$, we will denote a closed interval
\[
I_f:=
\left[ \inf \left\{x\in\bR\mid f(x-1)\ne f(x)-1\right\}
,
\sup \left\{x\in\bR\mid f(x+1)\ne f(x)+1\right\}
\right].\]\end{notation}
Here, we use the usual convention that $\inf \varnothing =\infty$ and $\sup\varnothing=-\infty$.
We see that each $f\in\Homeo_\bZ(\bR)$ satisfies $I_f=\varnothing$.
If $f\in\Homeo_c(\bR)$ then $I_f$ is the smallest compact connected set containing $\supp f$.

Let $f\in\Homeo_+(\bR)$.  We have that
\[
\begin{cases}
f(x+1)=f(x)+1,&\text{ for } x\ge \sup I_f;\\
f(x-1)=f(x)-1,&\text{ for } x\le \inf I_f.
\end{cases}\]
This observation justifies the following terminology.

\bd A homeomorphism $f\co\bR\to\bR$ is \emph{eventually $1$--periodic}  if $I_f$ is compact.\ed

For a compact interval $J\sse\bR$, we define
\begin{align*}
\Homeo_{\mathrm{ep}}(\bR)&:=\{f\in\Homeo_+(\bR)\mid f\text{ is eventually }1\text{--periodic}\};\\
\Homeo_{\mathrm{ep},J}(\bR)&:=\{f\in\Homeo_+(\bR)\mid I_f\sse J\}.
\end{align*}
It is obvious that $\Homeo_{\mathrm{ep}}(\bR)$ is a group with the function composition. However, $\Homeo_{\mathrm{ep},J}(\bR)$ is not a subgroup of $\Homeo_{\mathrm{ep}}(\bR)$ for any interval $J$ of positive length, since $T$ is an element of $\Homeo_{\mathrm{ep},J}(\bR)$ and $I_{f \circ T} = I_f - 1 $ for $f \in \Homeo_{\mathrm{ep}}(\bR)$.

For each $*\in\{c,\bZ,\mathrm{ep},J,(\mathrm{ep},J)\}$, we define \[
\Diff_*^{k[,\alpha]}(\bR)=\Diff_+^{k[,\alpha]}(\bR)\cap \Homeo_*(\bR).\]
Eventual periodicity generalizes the condition of having a compact support:
\[\Homeo_c(\bR)\cup\Homeo_\bZ(\bR)\sse\Homeo_{\mathrm{ep}}(\bR).\]
We also note that if $f\in\Homeo_{\mathrm{ep}}(\bR)$ then
the compact interval
\[
[\inf I_f-1,\sup I_f]\]
contains the set
\[
\supp [f^{-1},T^{-1}] = \overline{\{ x\in\bR\mid f(x+1)\ne f(x)+1\}}.\]
Actually, it is the smallest compact interval containing $\supp [f^{-1},T^{-1}]$.


For a map $f\co\bR\to\bR$, we write $\|f\|=\|f\|_\infty=\sup_{x\in\bR} |f(x)|$.
We also write
\[ \|f\|_k=\abss*{ \der{f}{k}},\quad \|f\|_{k,\alpha}=\brac*{\der{f}{k}}_\alpha,\] whenever these quantities are defined.

Let $f,g\in\Homeo_{\mathrm{ep}}(\bR)$.  We define their distance in the \emph{$C^0$--metric} as
\[
d_0(f,g)=\max(\abss{f-g}_0,\abss{f^{-1}-g^{-1}}_0).\]
If $k\ge1$, then the \emph{$C^{k[,\alpha]}$--metric} between $f,g\in \Diff_{\mathrm{ep}}^{k[,\alpha]}(\bR)$ is defined by
\begin{align*}
d_k(f,g)&=\sup_{0\le i\le k}\abss{f-g}_i;\\
d_{k,\alpha}(f,g) &=\max\left(\sup_{0\le i\le k}\abss{f-g}_i,   \sup_{1\le i\le k}\abss{f-g}_{i,\alpha}\right).
\end{align*}
The metric $d_{k[,\alpha]}$ determines the \emph{$C^{k[,\alpha]}$--topology} on $\Diff_{\mathrm{ep}}^{k[,\alpha]}(\bR)$.

For a group  $G$, we define the following maps $G\to G$:
\[
R_g(h)=h\circ g,\quad
L_g(h)=g\circ h,\quad
\Inv(h)=h^{-1}.\]

We note the following observations concerning the $C^{k[,\alpha]}$ topology, which we will prove in the course of this section:

\begin{prop}\label{prop:ck-group}
For each $k\in\bN$,  the following hold.
\be
\item\label{p:ck1}
The sets $\Diff^k_{\mathrm{ep}}(\bR)$ and $\Diff^{k,\alpha}_{\mathrm{ep}}(\bR)$ are topological groups in the $C^k$--topology. 
\item\label{p:ck2}
The right-multiplication map $R_g\co h\mapsto hg$ is $C^{k,\alpha}$--continuous in $\Diff^{k,\alpha}_{\mathrm{ep}}(\bR)$.
\item\label{p:ck3}
Let $J$ be a compact interval and let $*\in\{\mathrm{ep},\bZ,J\}$.
Then  $\Diff^{k[,\alpha]}_{*}(\bR)$ is connected in the $C^{k[,\alpha]}$--topology. 
\ee
\end{prop}

By Lemma~\ref{l:top-gp-gen}, we can deduce the following.
\begin{cor}\label{cor:ball-gen}
If $*\in\{\mathrm{ep},\bZ,J\}$, then every $C^{k[,\alpha]}$--neighborhood of the identity generates the group $\Diff^{k[,\alpha]}_*(\bR)$.
\end{cor}

\subsection{A remark on circles}\label{ss:rem-s1}
One can define the $C^{k[,\alpha]}$--metric on $\Diff_+^{k[,\alpha]}(S^1)$ in the same way as for eventually periodic diffeomorphisms
of the real line. 
The reader will note, for instance, that we have used the $C^1$--metric on $\Diff_+^1(S^1)$
in the course of the proof of Lemma~\ref{l:top-group}.

Lifting orientation preserving homeomorphisms of the circle to the real line, one obtains the \emph{universal central extension} 
\[1\to \form{T}\to\Homeo_\bZ(\bR)\to \Homeo_+(S^1)\to1.\]
The reader is directed to~\cite{Milnor1971book,Ghys2001,Calegari2007} for details.
If a map $h\co S^1\to\bR$ lifts to $\tilde h\co \bR\to\bR$, then the periodicity of $\tilde h$ implies
\[
[h]_\alpha =\sup_{x\ne y\in S^1}\frac{|hx-hy|}{\alpha(x-y)}
=\sup_{0<|x-y|\le 1/2} \frac{\abs*{\tilde hx-\tilde hy}}{\alpha(x-y)}
=\sup_{x\ne y} \frac{\abs*{\tilde hx-\tilde hy}}{\alpha(x-y)}=\brac*{\tilde h}_\alpha.\]

Suppose we have the following commutative diagram:
\begin{equation}\label{eq:lift}
\xymatrix{
\Diff^{k[,\alpha]}_\bZ(\bR)\ar[r]^{\tilde\Phi}\ar[d]^p&\Diff^{k[,\alpha]}_\bZ(\bR)\ar[d]^p\\
\Diff_+^{k[,\alpha]}(S^1)\ar[r]^\Phi &\Diff_+^{k[,\alpha]}(S^1)
}
\end{equation}
Here, the vertical maps $p$ are induced by the natural projection $\bR\to S^1=\bR/\bZ$.

If $p(\tilde f)=f$ and $p(\tilde g)=g$, then for each $i\ge1$ we have
\[
\abss{f}_{i[,\alpha]}=\abss*{\tilde f}_{i[,\alpha]},\quad
\abss{ f -  g}_{i[,\alpha]} = \abss*{\tilde f-\tilde g}_{i[,\alpha]}.\]
Moreover, if $\abss{\tilde f - \tilde g}\le 1/2$ then $\abss{\tilde f-\tilde g}=\abss{f-g}$.
It follows that the map $p$ is a $C^{k[,\alpha]}$--local-isometry.
In particular, 
 $\tilde\Phi$ is  $C^{k[,\alpha]}$--continuous, then so is $\Phi$.
By setting $\tilde\Phi$ as one of $L_g,R_g$ and $\operatorname{Inv}$, we obtain from Proposition~\ref{prop:ck-group} 
that the same statements as the first two parts also hold for $\Diff_+^{k[,\alpha]}(S^1)$. 

\subsection{$C^k$--continuity}
We prove the $C^k$--part of Proposition~\ref{prop:ck-group} in this subsection.
\begin{lem}\label{l:ep-gp}
The following hold.
\be
\item\label{p:epgp1}
The set  $\Diff^{k[,\alpha]}_{\mathrm{ep}}(\bR)$  is a group. 
\item\label{p:epgp2}
If $f\in\Diff^1_{\mathrm{ep}}(\bR)$, then 
$\inf f'>0$.
\item\label{p:epgp3}
For each $f\in\Diff^1_{\mathrm{ep}}(\bR)$ and $g\in\Homeo_{\mathrm{ep}}(\bR)$,
we have that 
\[
\abss{f^{-1}-g^{-1}}\le \abss{f^{-1}}_1\cdot \abss{f-g}.\]
\item\label{p:epgp4}
If $f\in\Diff^k_{\mathrm{ep}}(\bR)$
 for some $k\in\bN$, then 
 $\der{f}{i}$ is uniformly continuous for all $0\le i\le k$.
Moreover, for all $x\in \bR$ and $1\le i\le k$, 
we can find some $y_i$ with $|x-y_i|\le |I_f|/2+1$
such that $\der{f}{i}(y_i)=\der{\Id}{i}$.
\ee
\end{lem}
\bp 
Part (\ref{p:epgp1}) follows from
the fact that $\Homeo_{\mathrm{ep}}(\bR)$ is a group, together with Lemma~\ref{l:higher} and Lemma~\ref{l:inverse}.
For part (\ref{p:epgp2}), it suffices to note that $f$ is a diffeomorphism on any compact interval so that $f'$ is bounded
away from zero on any such interval, together with the fact that $f'$ is periodic outside a sufficiently large compact interval.
Part (\ref{p:epgp3}) can be seen from the estimate
\begin{align*}
\abss*{f^{-1}-g^{-1}}&=\abss*{f^{-1}-f^{-1}\circ(f\circ g^{-1})}\le\abss*{f^{-1}}_1\cdot \abss*{\Id-f\circ g^{-1}}\\
&\le\abss{f^{-1}}_1\cdot \abss*{g\circ g^{-1}-f\circ g^{-1}}=\abss{f^{-1}}_1\cdot \abss*{g-f}.\end{align*}

The uniform continuity in part (\ref{p:epgp4}) is obvious. For the remaining claim, recall $f\circ T(x)=T\circ f(x)$ for all $x$ outside $I_f$.
By the Mean Value Theorem, there exist some $y^+\in(\sup I_f,\sup I_f+1)$ and $y^-\in(\inf I_f-1,\inf I_f)$ such that 
\[f'((y^+ +\bZ)\cap(\sup I_f,\infty))=1=f'((y^- +\bZ)\cap(-\infty,\inf  I_f)).\]
Inductively, we can find $y_2,y_3,\ldots$ in $ [\sup I_f,\sup I_f+1]$ and also in $[\inf I_f-1,\inf I_f]$ such that $\der{f}{i}(y_i)=\der{\Id}{i}(y_i)$. 
\ep
Observe that part (\ref{p:epgp3}) of Lemma~\ref{l:ep-gp} implies that the $C^k$--topology is finer than the $C^0$--topology for all $k\in\bN$.
It will be convenient for us to use the notation
\[ M_k(f)=\sup_{1\le i\le k}\abss{f-\Id}_i=\max\left(\abss{f'-1},\sup_{2\le i\le k}\abss{f}_i\right).\]

\begin{lem}\label{l:ck-cont}
The following hold.
\be 
\item\label{p:rg}
Each right-multiplication map in $\Diff_{\mathrm{ep}}^k(\bR)$ is uniformly $C^k$--continuous. More precisely, there exists a constant $C=C(k)\in\bN$ 
such that for all $f_1,f_2,g\in \Diff_{\mathrm{ep}}^k(\bR)$ we have that
\[
d_k(R_g f_1,R_g f_2)\le C(1+ M_k(g))^k d_k(f_1,f_2).\]
\item\label{p:lg}
Each left-multiplication map and the inverse map in $ \Diff_{\mathrm{ep}}^k(\bR)$ are $C^k$--continuous. 
\ee
\end{lem}

\bp
Put $G= \Diff_{\mathrm{ep}}^k(\bR)$.
We let $f_1,f_2\in G $ be given, and set $a=f_1-f_2\co\bR\to\bR$. 

(\ref{p:rg})
We note that
\[
\abss{a\circ g}_k
\le
\abss*{
\left(\der{a}{k}\circ g\right)(g')^k}
+
\abss*{(a'\circ g)\cdot \der{g}{k}}
+
C_1
\sum_\gamma
\abss*{\left(a^{(i)}\circ g\right)\prod_{s=1}^i \der{g}{j_s}}
\]
for some $C_1=C_1(k)\ge0$. 
Here, $\gamma=({i,j_1,\ldots,j_i})$ varies over suitable indices.
Since $\abss{g}_i\le 1+\abss{g-\Id}_i$, we have
\[
\abss{a\circ g}_k\le
C'(1+ M_k(g))^k
\cdot
\sup_{1\le i\le k}\abss{a}_i\]
for some $C'\ge0$.
The (uniform) continuity of $R_g$ follows.

(\ref{p:lg})
Pick a concave modulus $\alpha$ such that  $g^{(k)}\in C^{\alpha}(\bR)$.
For $0\le i<k$, we have
\begin{align*}
\abss{\der{g}{i}\circ f_1-\der{g}{i}\circ f_2}&\le \abss{g}_{i+1}\cdot \abss{a},\\
\abss{\der{g}{k}\circ f_1-\der{g}{k}\circ f_2}&\le \abss{g}_{k,\alpha} 
\alpha\left(\abss{a}\right).
\end{align*}
For $m=1,2$, we can write
\[
\der{(g\circ f_m)}{k}=\left(\der{g}{k}\circ f_m\right)(f_m')^k+(g'\circ f_m)\cdot \der{f_m}{k}+\sum_\gamma C_\gamma  \left(g^{(i)}\circ f_m\right)\prod_{s=1}^i \der{f_m}{j_s}.\]
So, $d_k(L_gf_1,L_gf_2)$ is controlled by $d_k(f_1,f_2)=\sup_i \abss{a}_i$, whence the continuity of $L_g$ follows.

Fix $f_1\in G$. We prove that $d_k(f_1^{-1},f^{-1})$ can be made as small as we like by setting $d_k(f_1,f)$ to be sufficiently small.
For this, we put $\alpha_i(x) = \abss{f_1}_{i+1}\cdot x$ for $i<k$, 
and let $\alpha_k$ be a concave modulus for $\der{f_1}{k}$.
For compactness of notation, we put $F_1=f_1^{-1}$ and $F=f^{-1}$. 

We first require \[d_1(f_1,f )\le\frac12\inf f_1',\] so that $\inf f' \ge\frac12 \inf f_1'>0$.
It now suffices to note the following estimates:
\begin{align*} 
\abss{F_1'-F'}&\le \abss*{\frac{f_1'\circ F_1 - f'\circ F }{f_1'\circ F_1 \cdot f'\circ F } }\le \frac{2\abss{ f_1'\circ F_1- f'\circ F}}{\inf (f_1')^2},\\
 \abss{ f_1'\circ F_1- f'\circ F}&\le 
 \abss{f_1'\circ F_1- f_1'\circ F} +\abss{f_1-f}_1 \\
 &\le \alpha_1(\abss{F_1 - F})+\|a\|_1.\end{align*} 

Let $2\le r\le k$.
From Lemma~\ref{l:higher}, we have the identity
\[
0 = \der{(f \circ F)}{r} = (\der{f}{r}\circ F) (F')^{r+1} + (f' \circ F)\der{F}{r} + \sum_{\gamma \in A_r}{C_\gamma (\der{f}{i}\circ F)\prod_{t=1}^{i}\der{F}{j_t} }.
\]
Since $(f' \circ F)F' = (f \circ F)' = 1$, we get
\begin{equation}\der{F}{r}=-(\der{f}{r}\circ F)(F')^{r+1}- \left(\sum_{\gamma \in A_r}{C_\gamma (\der{f}{i}\circ F)\prod_{t=1}^{i}\der{F}{j_t} }\right)F',\end{equation}

where $C_*$ are integers and the set $A_r$ is the same as in Lemma~\ref{l:higher}. Thus, we obtained an estimate of the form

\begin{align*}
\abss*{F_1-F}_r\le
&\abss*{(\der{f_1}{r}\circ F_1)  (F_1')^{r+1}
-(\der{f}{r}\circ F)  (F')^{r+1}}
\\
&+\sum_{1<i<r,\alpha} C^*
\abss*{
\left(f_1^{(i)}\circ F_1\right)\prod_{s=0}^i \der{F_1}{j_s}
-\left(f^{(i)}\circ F\right)\prod_{s=0}^i \der{F}{j_s}
}.
\end{align*}
Inductively, we may assume $\abss{F_1-F}_j$ has been made as small as we desire for all $j<r$.
Then,
\[ \abss*{ \der{f_1}{i}\circ F_1- \der{f}{i}\circ F} \le \alpha_i(\abss{F_1-F}) + \abss{f_1-f}_i \] can also be made as small we we like
 small for $i\le r$.
It follows that $\Inv$ is continuous.
\ep

\bp[Proof of Proposition~\ref{prop:ck-group} (\ref{p:ck1})]
It only remains to show that the multiplication is $C^k$--continuous as a two-variable map $G\times G\to G$. Denote by $B(f;r)$ the $C^k$--ball of radius $r$ centered at $f\in G$. Let $f,g\in G$ be given, and let $\epsilon>0$ be arbitrary. Since $L_f$ is continuous
by Lemma~\ref{l:ck-cont} (\ref{p:lg}), we may suppose the existence of an $r_1\in(0,1)$ such that 
$f \circ B(g;r_1)\sse B(f \circ g;\epsilon)$.
For $f_1\in G$ near $f$ and $g_1 \in B(g;r_1)$, we have
\[
d_k(f_1\circ g_1, f\circ g_1)\le C(1+ M_k(g_1))^k d_k(f_1,f)\le C(2+ M_k(g))^k d_k(f_1,f).\]
So, we can find an $r_0>0$ such that 
\[B(f;r_0)\circ B(g;r_1)\sse B(f\circ g;2\epsilon).\qedhere\]
\ep

\subsection{$C^{k,\alpha}$--continuity}

\begin{lem}[Domination Lemma]\label{l:dom}
Let $J\sse\bR$ be a compact interval,
 let $f,g\in\Homeo_{\mathrm{ep},J}(\bR)$,
 and
 let $i$ be a nonnegative integer.
In the case when $i=0$, we further assume that  $f=g$ on $\partial J$.
Then for  $\ell=(|J|+2)$, there exists some $C=C(\alpha)$ such that
\[
\abss{f-g}_i\le C \ell\abss{f-g}_{i,\alpha},\quad
\abss{f-g}_i\le \ell \abss{f-g}_{i+1},\quad
 \abss{f-g}_{i,\alpha}\le C\ell\abss{f-g}_{i+1},\]
 whenever the relevant quantities are defined.
\end{lem}
\bp
Put $a=f-g$. Let us first assume $i\ge1$. 
For each $x\in \bR$ there exists $y\in \bR$ 
such that $|x-y|\le \ell$ and such that $\der{a}{i}(y)=0$, by part~(\ref{p:epgp4}) of Lemma~\ref{l:ep-gp}.
It follows that
\begin{align*}
\abs*{\der{a}{i}(x)} &=\abs*{\der{a}{i}(x)-\der{a}{i}(y)}\le \ell\alpha(1) \abss{a}_{i,\alpha}.\\
\abs*{\der{a}{i}(x)} &\le 
\int_y^x \abs*{\der{a}{i+1}}\le \ell \abss{a}_{i+1},\\
\abss{a}_{i,\alpha}
&=
\sup_{0<s-t\le\ell}
	\abs*
		{\frac
			{\der{a}{i}(s)-\der{a}{i}(t)}{\alpha(s-t)}
		}
\le
\sup_{0<s-t\le \ell}\frac{\abs{s-t}}{\alpha(s-t)}\abss{a}_{i+1}
\le \frac{\ell}{\alpha(\ell)}\abss{a}_{i+1}.
\end{align*}

For the case $i=0$,
The map $a=f-g$ satisfies that $a(\sup J+m)=0=a(\inf J-m)$ for all $m\in\bN\cup\{0\}$.
So, for each $x\in \bR$ there exists $y\in \bR$ 
such that $|x-y|\le \ell$ and such that $ay=0$.
The remainder of the argument is identical to the proof for the case $i>0$.
\ep

In particular, the metric $\abss{f-g}_{k[,\alpha]}$ determines the $C^{k[,\alpha]}$--topology on the group $\Diff_J^{k[,\alpha]}(\bR)$. Also, $\abss{f-g}=\abss{f-g}_0$ for $f,g\in\Diff_\bZ^{k[,\alpha]}(\bR)$ 
can be dominated by higher norms,
provided we require that $f(x)=g(x)$ for each $x\in\bZ$.


\begin{lem}[Derivation Formulae]\label{l:derivation}
For $f,g,a_1,\ldots,a_m\in C^\alpha(\bR)$, we have:
\be
\item
$[f\cdot g]_\alpha\le [f]_\alpha\abss{g}+\abss{f}[g]_\alpha$.
\item
$\brac{a_1  a_2\cdots a_m}_\alpha\le \sup_i \abss{a_i}^{m-1}\sum_i[a_i]_\alpha$.
\item\label{p:fgalpha}
Additionally, 
if $g\in \Diff_+^1(\bR)$, then
$[f\circ g]_\alpha\le [f]_\alpha \max(\abss{g}_1,1)\le \brac{f}_\alpha(1+M_1(g))$.
\ee
\end{lem}
\bp
We may assume the right-hand side of each item in the statement of the lemma is finite. We use the notation
\[
\eval{f}^x_y:=f(x)-f(y).\]
For part (1), note that
\[
\frac{\eval{f\cdot g}^x_y}{\alpha(x-y)}
=
\frac{\eval{f}^x_y}{\alpha(x-y)} g(x)
+
f(y) \frac{\eval{g}^x_y}{\alpha(x-y)}
.\]
Part (2) follows by a straightforward induction.
For part (\ref{p:fgalpha}), we note $|gx-gy|\le\abss{g}_1 |x-y|$
and 
\[
\frac{\abs{\eval{f\circ g}^x_y}}{\alpha(x-y)}
=
\frac{\abs{\eval{f}^{gx}_{gy}}}{\alpha(gx-gy)}\cdot\frac{\alpha(gx-gy)}{\alpha(x-y)}
\le [f]_\alpha \max(\abss{g}_1,1).\qedhere\]
\ep


\bp[Proof of Proposition~\ref{prop:ck-group} (\ref{p:ck2})]
We continue from the proof of Lemma~\ref{l:ck-cont}.
In particular, we let $f_1,f_2\in \Diff^{k,\alpha}_{\mathrm{ep}}(\bR)$ and put $a=f_1-f_2$. 
For some $K=K(k,g)$ we have that
\[
d_k(f_1\circ g, f_2\circ g)\le K d_k(f_1,f_2).\]
Also, there exists some $C=C(k)\ge0$ such that
\[
\abss{a\circ g}_{k,\alpha}
\le
\brac*{\left(\der{a}{k}\circ g\right)(g')^k}_\alpha
+\brac*{(a'\circ g) \der{g}{k}}_\alpha
+C\sum_{\gamma} \brac*{\left(a^{(i)}\circ g\right)\prod_{s=1}^i \der{g}{j_s}}_\alpha.\]
By the Derivation Formulae (Lemma~\ref{l:derivation}), we estimate each term as 
\begin{align*}
\brac*{\left(a^{(i)}\circ g\right)\prod_{s=1}^i \der{g}{j_s}
}_\alpha&
\le \brac*{\der{a}{i}\circ g}_\alpha\prod_{s=1}^i\abss{g}_{j_s}
+\abss{a}_i\brac*{\prod_{s=1}^i\der{g}{j_s}}_\alpha.\\
\brac*{\der{a}{i}\circ g}_\alpha&\le \abss{a}_{i,\alpha}\max(\abss{g}_1,1).
\end{align*}
So, we have $K = K(k,g,\alpha)$ such that
\[\abss{(f_1-f_2)\circ g}_{k,\alpha}=\abss{a\circ g}_{k,\alpha}\le K d_{k,\alpha}(f_1,f_2).\]
By the Domination Lemma, we obtain that 
$d_{k,\alpha}(f_1\circ g, f_2\circ g)$ can be made arbitrarily small.
\ep

\begin{exmp}
The map $L_g$ is not necessarily $C^{k,\alpha}$--continuous even at $f=\Id$. For any concave modulus $\alpha$, pick $f_1,g \in \Diff_c^{k,\alpha}(\bR)$ such that $f_1(x) = x+1$ and $\der{g}{k}(x) = \alpha(|x|)$ near $x=0$, and set $f_t (x) = tf_1(x) + (1-t)\Id(x)$ for $0 \le t \le 1$. Then, $f_t$ converges to $\Id$ as $t \rightarrow 0$, but 
\[
\der{L_g(f_t)}{k} - \der{L_g(\Id)}{k} = \alpha(|x+t|) - \alpha(|x|),
\]
\[
\brac*{\der{L_g(f_t)}{k} - \der{L_g(\Id)}{k}}_\alpha \ge \frac{1}{\alpha(0-(-t))} (\alpha(|t|) + \alpha(|-t|)) = 2,
\]
which implies that $L_g(f_t)$ does not converges to $L_g(\Id)$.
\end{exmp}

Part (\ref{p:ck3}) of Proposition~\ref{prop:ck-group} is trivially is implied by the following:

\begin{lem}\label{l:top-group-komega} 
Let $*$ be as in Proposition~\ref{prop:ck-group}. Then the following metric spaces are contractible:
\[\left(\Diff^{k}_*(\bR), d_{k}\right),\quad \left(\Diff^{k,\alpha}_*(\bR), d_{k,\alpha}\right).\]\end{lem}
\bp
Set $G = \Diff^{k[,\alpha]}_*(\bR)$.
We claim $g_t(x):=H(t,g)(x)=tg(x)+(1-t)x$ defines a  deformation retract
 \[H\co [0,1]\times  G\to  G,\quad H(0,G)=\{\Id\},\quad H(1,g)=g.\]
To see this, first observe from $\partial g_t(x)/\partial x=tg'(x)+(1-t)>0$ that $g_t\in  G$.
Since the following hold for each $i=0,1,\ldots,k$, we conclude that  $H$ is continuous.
\begin{align*}
&\abss{g_s-g_t}_{i[,\alpha]}\le |s-t|\cdot 
\abss{g}_{i[,\alpha]},\\
&\abss{g_t - h_t}_{i[,\alpha]}\le 
t \abss{g-h}_{i[,\alpha]}.\qedhere
\end{align*}
\ep

\begin{rem}
For a given $g\in \Homeo_*(\bR)$ and $s\in [0,1]$ we can find
 a concave modulus $\alpha$ such that $[g_s^{-1}]_\alpha<\infty$; see 
 Theorem~\ref{thm:decompose}.
 Then, we have
\[
\abss{g_s^{-1}-h_t^{-1}}
=\abss{g_s^{-1}\circ h_t-g_s^{-1}\circ g_s}
\le \brac{g_s^{-1}}_\alpha\cdot \alpha(\abss{g_s- h_t}).\]
Combining with the proof of Lemma~\ref{l:top-group-komega},
we see that $\Homeo_*(\bR)$ is also  contractible  (in the $C^0$--topology).
\end{rem}

\subsection{Compositions and norms}\label{subsec:comp-and-norms}
Let $J$ be a compact interval.
Let $*\in\{c,\bZ,\mathrm{ep},J,(\mathrm{ep},J)\}$ and $\delta>0$, where $J\sse\bR$. For $k\in\bN$, it will be convenient to write
\begin{align*}
V_*^0(\delta)&= \{ f\in \Homeo_*(\bR) \co \|f-\Id\|<\delta\};\\
V_*^{k[,\alpha]}(\delta)&=\{ f\in \Diff_*^{k[,\alpha]}(\bR)\co \abss{f-\Id}_{k[,\alpha]}<\delta\};\\
\bar V_*^{k[,\alpha]}(\delta)&=\{ f\in \Diff_*^{k[,\alpha]}(\bR)\co \abss{f-\Id}_{k[,\alpha]}\le \delta\}.
\end{align*}

The goal of this subsection is the following $C^{k,\alpha}$--estimate:
\begin{prop}\label{prop:muk-compose}
Let $k\in\bN$, and let $J\sse \bR$ be a compact interval. 
Then there exist constants $\epsilon =\epsilon (k,\alpha,|J|)$ and 
$C=C(k,\alpha,|J|)$
such that whenever $f,g\in V_{\mathrm{ep},J}^{k,\alpha}(\epsilon)$,
we have that 
 \[\abss{f\circ g}_{k,\alpha}\le \abss{f}_{k,\alpha}+\abss{ g}_{k,\alpha}
+C\abss{f}_{k,\alpha} \abss{ g}_{k,\alpha}.\]
\end{prop}




\begin{rem}\label{r:domination}
By the Domination Lemma, we see that the metric $\abss{f-g}_{k[,\alpha]}$ determines the $C^{k[,\alpha]}$--topology on the group $\Diff_J^{k[,\alpha]}(\bR)$. Also, $\abss{f-g}=\abss{f-g}_0$ for $f,g\in\Diff_\bZ^{k[,\alpha]}(\bR)$ 
can be dominated by higher norms
say, if we further require that $f(x)=g(x)$ for each $x\in\bZ$.
\end{rem}

\bp[Proof of Proposition~\ref{prop:muk-compose}]
The Domination Lemma implies that $M_k(g)\le 1$ if $\epsilon $ is chosen to be sufficiently small.
Here we recall the notation \[ M_k(g)=\sup_{1\le i\le k}\abss{g-\Id}_i=\max\left(\abss{g'-1},\sup_{2\le i\le k}\abss{g}_i\right).\]

Moreover,  we have some constant $C'$ such that 
\[
M_{k,\alpha}(f)  M_{k,\alpha}(g)\le
C'\abss{f}_{k,\alpha} \abss{ g}_{k,\alpha}.\]

Thus, it suffices to show the following claim:
\begin{claim*}
If $\epsilon$ and $C$ are suitably chosen a priori, then we have
 \[\abss{f\circ g}_{k,\alpha}\le \abss{f}_{k,\alpha}+\abss{ g}_{k,\alpha}
+C M_{k,\alpha}(f)  M_{k,\alpha}(g).\]
\end{claim*}
To prove the claim, we first note from $\abss{g'-1} \le 1$ that
\[
\abss{g}_1^{k+1}\le (1+\abss{g'-1} )^{k+1}\le 1+D_0\abss{g'-1}\]
for some $D_0=D_0(k)>0$.

The case $k=1$ follows from the Derivation Formulae (Lemma~\ref{l:derivation}) and from the computation 
\begin{align*}
\abss{f\circ g}_{1,\alpha}&=\brac*{(f'\circ g)g'}_{\alpha}
\le
\brac*{f'\circ g}_{\alpha}\cdot \abss{g}_1+\abss{f}_1\cdot \brac*{g'}_{\alpha}
\\
&\le  \abss{f}_{1,\alpha}\abss{g}_1^2 + (1+\abss{f'-1})\abss{g}_{1,\alpha}\\
&\le  \abss{f}_{1,\alpha}+\abss{g}_{1,\alpha}+ C_1  M_{1,\alpha}(f)M_{1,\alpha}(g).
\end{align*}
for some $C_1>0$.

Assume $k\ge2$. There is a constant $D=D(k)$ such that
\[
\abss{f\circ g}_{k,\alpha}
\le
\brac*{\left(\der{f}{k}\circ g\right)(g')^k}_{\alpha}
+
\brac*{(f'\circ g)\der{g}{k}}_{\alpha}+
D\sum_{i,j_s}  \brac*{( \der{f}{i}\circ g)\der{g}{j_1}\cdots\der{g}{j_i}}_{\alpha}.\]

We estimate each of the above three terms separately.
For the first term, we note  there exist constants $D_1=D_1(k),D_2=D_2(k)$ such that
\begin{align*}
\brac*{\der{f}{k}\circ g}_\alpha\abss*{(g')^k} &\le \abss{f}_{k,\alpha}(1+D_0\abss{g'-1} )
\le \abss{f}_{k,\alpha} + D_1 M_{k,\alpha}(f) M_{k,\alpha}(g),\\
\abss*{\der{f}{k}\circ g}\brac*{(g')^k}_\alpha
&\le D_2 \abss{f}_k \abss{g}_{1,\alpha}
\le D_2 M_{k,\alpha}(f) M_{k,\alpha}(g).
\end{align*}

For the second term, we find a constant $D_3$ such that
\begin{align*}
\brac*{f' \circ g}_\alpha\abss*{\der{g}{k}} &\le \abss{f}_{1,\alpha}(1+\abss{g'-1} )\abss{g}_k
\le
D_3 M_{k,\alpha}(f) M_{k,\alpha}(g).
\\
\abss*{f'\circ g}\brac*{\der{g}{k}}_\alpha
&\le (1+\abss{f'-1})\abss{ g}_{k,\alpha}\le \abss{ g}_{k,\alpha}+M_{k,\alpha}(f) M_{k,\alpha}(g).
\end{align*}

Let $(i,j_1,\ldots,j_i)$ be as in the third term, so that $1<i<k$ and $j_i>1$. Then
\begin{align*}
\brac*{\der{f}{i}\circ g}_\alpha&\le \abss{f}_{i,\alpha}(1+\abss{g'-1} )\le 2M_{k,\alpha}(f),\\
\abss*{\prod_{s=1}^i \der{g}{j_s}}&\le (1+ M_{k-1}(g))^{k-2}\cdot \abss{g}_{j_i}
\le 2^k  M_{k,\alpha}(g),\\
\abss{\der{f}{i}\circ g}&\le M_{k,\alpha}(f),\\
\brac*{\prod_{s=1}^i \der{g}{j_s}}_\alpha &\le
 \brac*{\der{g}{j_i}}_\alpha\cdot 2^{k-2}+ \abss{g}_{j_i}
\brac*{\prod_{s=1}^{i-1} \der{g}{j_s}}_\alpha
\le 2^{k-2}
 M_{k,\alpha}(g)
+
k\cdot 2^k M_k(g) M_{k,\alpha}(g).
\end{align*}
The claim follows.
\ep

\subsection{The fixed point property}
We say a topological space $X$ has the \emph{fixed point property} if every continuous map $T\co X\to X$ has a fixed point. 
The ingredient of the proof of Mather's Theorems we establish in this subsection is the following:
\begin{prop}\label{prop:fixed}
For a compact interval $J$ and for all sufficiently small $\epsilon>0$, 
the set
\[ \bar  V_J^{k,\alpha}(\epsilon).\]
has the fixed point property in the $C^k$--topology.
\end{prop}

We will require the following general result functional analysis:

\begin{thm}[{Schauder--Tychonoff Fixed Point Theorem~\cite[Chapter V]{DSbook1}}]\label{thm:schauder-fixed}
Every nonempty compact convex subset of a locally convex topological vector space (LCTVS)
has the fixed point property.\end{thm}

A precise definition of a LCTVS is not relevant for our purposes, so we will avoid any extended digression on this subject.
At the end of this section, we prove Theorem~\ref{thm:schauder-fixed}
for a reasonably large subclass of LCTVS, namely  \emph{Fr\'echet spaces}. 
For our purposes, namely proving Theorem~\ref{thm:mather-1d}, it would suffice to only consider Banach spaces.
However, in order to fully establish the $C^{\infty}$ case, 
one needs to apply Theorem~\ref{thm:schauder-fixed} to the Fr\'echet space  $C^\infty_J(\bR)$; see Remark~\ref{r:cinfty} and~\cite{Epstein1984CMH}.

For a compact interval $J$, we let  \[C_J^{k[,\alpha]}(\bR):=\{ f\in C^{k[,\alpha]}(\bR)\mid \supp f \sse J\}.\] 
The vector space $C_J^{k[,\alpha]}(\bR)$ is equipped with the \emph{$C^{k[,\alpha]}$--norm} $\abss{\cdot}_{k[,\alpha]}$.

\begin{lem}\label{l:lip-met}
Let $k\in\bN\cup\{0\}$. For $K\ge |J|+\alpha(|J|)+\frac{|J|}{\alpha(|J|)}$, the following hold:
\be
\item If $f\in C_J^{k,\alpha}(\bR)$, then $\abss{f}_k\le K \abss{f}_{k,\alpha}$.
\item
If $f\in C_J^{k+1}(\bR)$, then  $\abss{f}_k\le K \abss{f}_{k+1}$ and $  \abss{f}_{k,\alpha}\le K \abss{f}_{k+1}$.
\ee
\end{lem}
\bp
We may assume $k=0$. For $x \in J=[0,|J|]$, we have that
\begin{align*} 
\abs{fx}&=\abs{fx-f0}\le    \alpha(x) [f]_\alpha \le  \alpha(|J|)[f]_\alpha,\\
\abs{fx}&\le \abs*{\int_0^x f'}\le |J|\cdot \abss{f'},\\
\sup_{0<x-y<\infty}\frac{\abs{fx-fy}}{\alpha(x-y)}
&\le
\abss{f'}\cdot\sup_{0<x-y\le |J|}\frac{\abs{x-y}}{\alpha(x-y)}\le \frac{|J|}{\alpha(|J|)}\abss{f'}.\qedhere\end{align*}
\ep

The following two results are standard and can be found in any number of standard textbooks, so we omit any further discussion.

\begin{thm}[Arzel\`a--Ascoli]\label{thm:aa}
Every uniformly bounded, equicontinuous family of continuous functions on a compact metric space is relatively compact.
\end{thm}

\begin{lem}\label{l:c-banach0}
For $k\ge0$, the normed space $(C_J^k(\bR),\abss{\cdot}_{k})$ is Banach.
\end{lem}

Proposition~\ref{prop:fixed} follows from Schauder--Tychonoff Fixed Point Theorem and from Lemma~\ref{l:ball-cpt-d} below.
\begin{lem}\label{l:ball-cpt-d}
For real numbers $r\in(0,\infty), \epsilon\in(0,1)$ and for a compact interval $J\sse\bR$, the set
\[\{g\in \Diff_J^{k,\alpha}(\bR)\co  \abss{g}_{k,\alpha}\le r\text{ and }\abss{g'-1}\le\epsilon\}\]
is $C^k$--isometric to a compact convex subset of the Banach space $C_J^{k,\alpha}(\bR)$.
\end{lem}

\bp
We let $D$ be the given set, and put
\[\{f\in C_J^{k,\alpha}(\bR)\co  \abss{f}_{k,\alpha}\le r\text{ and }\abss{f'}\le\epsilon\}.\]
The map $\operatorname{VF}\co\DD(J;1)\to \CC(J;1)$ defined by 
$g\mapsto g-\Id$ is a $C^k$--isometric embedding,
which maps $D$ to $C$.
This shows that $\operatorname{VF}(D)=C$.
Since $[\cdot]_{k,\alpha}$ and $\abss{\cdot}_1$ are norms, we see that $C$ is convex.

It now suffices to show that $C$ is sequentially compact, as $C_J^k(\bR)$ is a metric space.
Pick a sequence $\{f_n\}\sse C$, for which we will find a convergent subsequence.
We may put $J=[0,\ell]$.

From Lemma~\ref{l:lip-met}, there is a constant $K=K(k,\alpha,\ell)>0$ such that \[\abss{f_n}_i\le Kr,\quad \abss{f_n}_{i,\alpha}\le Kr\]
for all $0\le i\le k$.
In particular, the class
$\FF_i:=\{\der{f_n}{i}\}$ is uniformly bounded and equicontinuous.
By the Arzel\`a--Ascoli Theorem,
we may assume that $\der{f_n}{i}$ uniformly converges to some $h_i\in C(J)$.
After enlarging $J$ if necessary, we may assume a priori that $f_n=0$ on some $\delta$--neighborhood of $\partial J$ for all $n$.
In particular, $h_i$ extends by $0$ outside $J$.

For each $i<k$ and $x\in J$, we note that
\[
h_i(x) =\lim\der{f_n}{i}(x) =\lim\int_0^x \der{f_n}{i+1} =\int_0^x h_{i+1},\] by the Dominated Convergence Theorem.
Hence, $h_i'=h_{i+1}$. It follows that $h_0\in C_J^k(\bR)$ and $f_n\to h_0$ in the $C^k$--metric.

Write $f=h_0$.
For $x,y\in J$ satisfying $0<x-y\le 1$, we have that
\[
\abs*{\frac{\der{f}{k}(x)-\der{f}{k}(y)}{\alpha(x-y)}}=\lim_n \abs*{\frac{\der{f_n}{k}(x)-\der{f_n}{k}(y)}{\alpha(x-y)}} \le r.\]
Thus, we have $\abss{f}_{k,\alpha}\le r$ and $f\in C_J^{k,\alpha}(\bR)$.
Furthermore, 
\[
|f'(x)| =\lim |f_n'(x) | \le \epsilon.\]
Hence, $f\in C$. This implies that $B$ is $C^k$--compact.
\ep

By Lemma~\ref{l:lip-met}, the condition $\abss{f'-1}\le\epsilon$ is redundant for a small $r>0$.

\subsection{Proof of the Schauder--Tychonoff Fixed Point Theorem}
Let us prove Theorem~\ref{thm:schauder-fixed} for Fr\'echet spaces.
Recall a \emph{Fr\'echet space} is a topological vector space $X$ 
such that for some family of seminorms $\{\abss{\cdot}_k\}_{k\in\bN_0}$,
 the topology of $X$ is determined by some complete metric $d$
given by \[d(f,g):=\sum_{k\ge0} 2^{-k}\frac{1}{1+1/\abss{f-g}_k}.\]

An important example of a (non-Banach) Frech\'et space is
\[ C_J^\infty(\bR):=\{f\co \bR\to\bR\mid \supp f\sse J\text{ and }f\text{ is }C^\infty\}.\]

In order to prove Theorem~\ref{thm:schauder-fixed}, we let
$X$ be a Frech\'et space whose topology is given by the family of seminorms $\FF=\{\abss{\cdot}_k\}_{k\ge1}$. We let $d$ be the induced metric from $\FF$ as above.
For each $n\in\bN$, we define \[d_k(x,y):=\sup_{i\le k} \abss{x-y}_i.\]

\begin{lem}\label{l:dkdense}
If $K\sse X$ is compact, then for each $k\in\bN$ and for each $\epsilon>0$
there exists a finite set $A$ such that 
\[K\sse\{ x\in X\mid d_k(x,A)<\epsilon\}.\]
\end{lem}
\bp
We choose $\delta>0$ so small that 
\[
0< \frac1{\frac1{2^k\delta}-1}<\epsilon.\]
Since $K$ is totally bounded with respect to $d$, there exists a finite $\delta$--net $A\sse K$.
For each $x\in K$, there exists $a\in A$ such that 
\[\sum_{i\le k} 2^{-i}\frac1{1+\frac1{\abss{a-x}_i}}\le d(a,x)<\delta.\]
For each $i\le k$, we have that
\[
\abss{a-x}_i\le \frac1{\frac1{2^i\delta}-1}\le \frac1{\frac1{2^k\delta}-1}<\epsilon.\]
Thus, we obtain $d_k(a,x)<\epsilon$.
\ep

We note that $d_k$ is \emph{convex}. Namely, if $x,x',y,y'\in X$ and $t\in[0,1]$ then
\begin{align*}
&d_k(tx+(1-t)y,tx'+(1-t)y')=\sup_i \abss{t(x-x')+(1-t)(y-y')}_i\\
&\le \sup_i t\abss{x-x'}_i+(1-t)\abss{y-y'}_i\le
t d_k(x,x')+(1-t) d_k(y,y').\end{align*}

For a finite subset $A$ of $X$, denote by $\hull(A)$ a convex hull of $A$. That is,
\[
\hull(A) = \{ \Sigma_{i=1}^m t_i a_i \mid \Sigma_{i=1}^m t_i = 1, t_i \ge 0, a_i \in A \}.\]

\begin{lem}[Schauder Projection Lemma]
Let $k\in\bN$ and $\epsilon>0$.
Suppose $A$ is a finite subset of a convex set $K\sse X$ such that 
\[K\sse\{ x\in X\mid d_k(x,A)<\epsilon\}.\]
Then there exists a continuous map
\[
\pi\co K\to \hull(A)\]
such that $\sup_K d_k(\pi x,x)\le\epsilon$.\end{lem}
\bp
We write $A=\{a_1,\ldots,a_m\}\sse K$.
Define $\psi_i\co K\to [0,\epsilon]$ by \[\psi_i(v)=\max(0,\epsilon-d_k(v,a_i)).\]
We have that $\psi=\sum_i\psi_i\co K\to(0,\infty)$. We define a continuous map $\pi\co K\to\hull(A)$ by  $\pi = \sum_i(\psi_i/\psi)a_i$. 
Note that $\psi_i(v)\ne0$ only if $d_k(a_i,v)< \epsilon$.
From the convexity of $d_k$, we have
 \[d_k(\pi(v),v)=d_k\left(\sum_i \frac{\psi_i(v)}{\psi(v)}a_i , \sum_i \frac{\psi_i(v)}{\psi(v)}v\right)
\le \sum_i \frac{\psi_i(v)}{\psi(v)}d_k(a_i, v)\le\epsilon\sum_i \frac{\psi_i(v)}{\psi(v)}=\epsilon.\qedhere\]
\ep

\bp[Proof of Theorem~\ref{thm:schauder-fixed}]
For each $k\in\bN$, we pick $A_k\sse K$ and $\pi_k\co K\to\hull(A_k)$ as in Lemma~\ref{l:dkdense} and the Schauder Projection Lemma,
after setting $\epsilon=1/k$.
Consider $T_k\co\hull(A_k)\to\hull(A_k)$ defined as the composition
\[
\xymatrix{
\hull(A_k)\ar@{^(->}[r]& K\ar[r]^T & K\ar[r]^<<<<{\pi_k} &\hull(A_k).}\]
Note that there exists a linear isomorphism between $\form{A_k}\le X$ and a finite dimensional Euclidean space, which is a homeomorphism. We remark that this fact holds in a more generalized setting; see~\cite{Rudin1991book}, for instance.

Thus, we have that, $\hull(A_k)$ is homeomorphic to a Euclidean ball.
The Brouwer Fixed Point Theorem implies that $\pi_k(T v_k)=v_k$ for some $v_k\in\hull(A_k)\sse K$.

Since $K$ is (sequentially) compact, we have $v_k\to v\in K$  after taking a subsequence.
For each $m\ge k\ge1$, we compute
\begin{align*}
&\abss{T(v)-v}_k\le \abss{T(v)-T(v_m)}_k + \abss{T(v_m)-\pi_m T(v_m)}_k
+\abss{v_m-v}_k,\\
&\abss{T(v_m)-\pi_m T(v_m)}_k\le d_m(T(v_m),\pi_m T(v_m))\le 1/m,\\
&\lim_{m\to\infty} d_k(v,v_m)=0=\lim_{m\to\infty}d_k(Tv,Tv_m).
\end{align*}
This implies that $\abss{T(v)-v}_k=0$ for all $k$, thus furnishing the desired fixed point.\ep

\section{The norm reduction operator}\label{s:norm}
Throughout the remainder of this article we consider a universal constant $\epsilon_0\in(0,1/100)$, which,
when the need arises, will be further shrunk.
We continue to let $k$ be a positive integer, and let $\alpha$ be a concave modulus.
\subsection{The main construction}
The heart of Mather's argument is the construction of
a certain operator $\Psi$ that reduces the $C^{k,\alpha}$--norm of a given diffeomorphism
without changing the first homology class. 

Let $J\sse\bR$ be a compact interval. For $*\in\{c,\bZ,\mathrm{ep},J\}$,
 we have defined
\begin{align*}
V_*^{k[,\alpha]}(\delta)&=\{ f\in \Diff_*^{k[,\alpha]}(\bR)\co \abss{f-\Id}_{k[,\alpha]}<\delta\};\\
V_{\mathrm{ep},J}^{k[,\alpha]}(\delta)&
=\{f\in  V_{\mathrm{ep}}^{k[,\alpha]}(\delta)\mid
\supp [f^{-1},T^{-1}]\sse J\}.
\end{align*}
We also let $\bar V_*$ be defined from $V_*$ as above by replacing ``$<\delta$'' by ``$\le\delta$''.
\begin{thm}\label{t:mather-inv}
Let $A\in\bN$, and set
\begin{itemize}
\item $D=[-2,2]$ and $E=[-2A,2A]$ for $k\ge2$;
\item $D=[-2A,2A]$ and $E=[-2,2]$ for $k=1$.
\end{itemize}
Then there exist positive numbers $C=C(k,\alpha)$, $\delta=\delta(k,A,\alpha)$,
and  a $C^k$--continuous map
\[\Psi \co V_E^{k,\alpha}(\delta)\to\Diff_D^{k,\alpha}(\bR)\]
such that 
whenever $g\in V_E^{k,\alpha}(\delta) $ the following hold.
\be
\item $\abss{\Psi (g)}_{k,\alpha}\le C|E|/|D|\cdot \abss{g}_{k,\alpha}$;
\item $[\Psi(g)]=[g]$ in $H_1(\Diff_c^{k,\alpha}(\bR))$.
\ee
\end{thm}
\begin{rem}
We will actually show that 
\[
\tau  \Psi (g)=\lambda\tau g\lambda^{-1}\]
for some  $\tau\in\Diff_c^\infty(\bR)$ and $\lambda=\lambda(g)\in \Diff_c^{k,\alpha}(\bR)$. The conclusion (2) will then follow.
\end{rem}

\begin{cor}\label{c:diff-perf}
For each tame pair $(k-1,\alpha)$, the group $\Diff_c^{k,\alpha}(\bR)$ is perfect.
\end{cor}

\bp[Proof of Corollary~\ref{c:diff-perf}, assuming Theorem~\ref{t:mather-inv}]
Consider $A\in\bN$, whose value will be determined later. 
Let $\Psi,C,\delta, D, E$ be as in Theorem~\ref{t:mather-inv}, and let $\epsilon>0$ be chosen to be sufficiently small.
By Lemma~\ref{l:top-gp-gen} and Proposition~\ref{prop:ck-group} it suffices to
establish that an arbitrary element $f \in V_D^{k,\alpha}(\epsilon)$ is trivial in the first homology group of $\Diff_c^{k,\alpha}(\bR)$.

We have seen in Lemma~\ref{l:ball-cpt-d} that the following set is $C^k$--compact:
 \[ K =  \bar V_D^{k,\alpha}(\epsilon).\]
Let $u\in K$ be given. By shrinking $\epsilon$ if necessary, we have from Proposition~\ref{prop:muk-compose} that 
$\abss{fu}_{k,\alpha}\le 3\epsilon$.
We can find a diffeomorphism $Q\in \Diff_c^{\infty}(\bR)$ so that
$g:=Q fu Q^{-1}$ satisfies 
\[ g(x)=  \frac{|E|}{|D|}  \cdot f u \left(\frac{|D|}{|E|}\cdot x\right)\in\Diff_E^{k,\alpha}(\bR).\]

To construct such a $Q$,
we start with the diffeomorphism of $\bR$ given by scaling by the factor $|E|/|D|$,
and declare $Q$ to coincide with this scaling on some sufficiently large compact
interval containing $D$. To make $Q$ compactly supported, we simply choose a $C^{\infty}$ cutoff, so that outside a larger compact interval, $Q$ becomes
the identity.

We claim that if $A$ is chosen to be sufficiently large to begin,
then
 $\abss{\Psi g}_{k,\alpha}<\epsilon$.
 We first assume $k\ge2$.
Since  $(k-1,\alpha)$ is a tame pair and $|E|/|D|=A$, we can find $A\gg0$ such that
\[
\abss{\Psi g}_{k,\alpha}
\le C A\abss{g}_{k,\alpha}
\le CA \cdot A^{1-k}\left( \sup_{s>0}\frac{\alpha(s/A)}{\alpha(s)}\right)\cdot \abss{fu}_{k,\alpha}<\epsilon,\]
as follows from Theorem~\ref{t:mather-inv}. Here, it is crucial  that $C$ does not depend on the choice of $A$.
The case $k=1$ is very similar. Indeed, set $A\gg0$. From the hypothesis that $\alpha$ is sup-tame, we see that
\[\abss{\Psi g}_{k,\alpha}
\le \frac{C}{A}\cdot \abss{g}_{k,\alpha}
\le C \left( \sup_{s>0}\frac{\alpha(As)}{A \alpha(s)}\right)\cdot \abss{fu}_{k,\alpha}<\epsilon.\]

Thus, we have a $C^k$--continuous map $\Theta\co K\to K$ defined by composition:
\[
u\mapsto fu \mapsto g = Q f uQ^{-1} \mapsto \Psi g.\]
By Schauder--Tychonoff Fixed Point Theorem (Theorem~\ref{thm:schauder-fixed}), the map $\Theta$ fixes some $u_0\in K$.
It follows that
\[[u_0] =[\Theta u_0]= [fu_0]\in H_1(\Diff_c^{k,\alpha}(\bR);\bZ).\] This implies that $[f]=0$.\ep

\bp[Proof of Theorem~\ref{thm:mather-1d}, assuming Theorem~\ref{t:mather-inv}]
Let $G=\Diff_c^{k,\alpha}(\bR)$ or  $G=\Diff_+^{k,\alpha}(S^1)$.
Applying fragmentation, each diffeomorphism $g\in \Diff_+^{k,\alpha}(S^1)$
can be written as a product of diffeomorphisms supported on open subintervals. 
It follows from Lemma~\ref{c:diff-perf} that $G$ is perfect.
As $[G,G]$ is simple (Section~\ref{s:higman}), we see that $G$ is simple.
Corollary~\ref{c:decompose} similarly implies that the groups $G=\Diff_c^k(\bR)$ or  $G=\Diff_+^k(S^1)$ are simple.

Since $\beta(t) = t^s$ is both a sup-tame and sub-tame concave modulus for $s \in (0,1)$,
we see that both $\Diff_c^s(\bR)$ and $\Diff_+^s(S^1)$ are perfect for every $s>1$.
Thus, for $r\ge1$, the groups
\[\bigcup_{s>r} \Diff_c^s(\bR)\] and \[\bigcup_{s>r}\Diff_+^s(S^1)\] are perfect, and their commutators are simple by the
Epstein-Ling Theorem (Theorem~\ref{thm:epstein-ling}). This shows that both of  these groups are simple.
Corollary~\ref{c:decompose} now implies that the groups
\[G=\bigcap_{s<r} \Diff_c^s(\bR)\] and \[G=\bigcap_{s<r}\Diff_+^s(S^1)\] are simple for $r>3$.
\ep

\begin{rem}\label{r:cinfty}
Epstein followed the same approach in the proof of the perfectness of $\Diff_c^\infty(\bR)$.
Through a more detailed analysis of  the norms $\abss{\Psi(g)}_k$,
he proved that the same map $\Theta$ in the proof of Corollary~\ref{c:diff-perf} is a continuous (in the $C^\infty$--topology) self-map of the set
\[
L=\{g\co  \Diff_D^\infty(\bR) \co \abss{g}_i\le \epsilon_i\text{ for all }i\ge 3\}\]
for a suitable choice of $\{\epsilon_i\}$.
Since $L$ is naturally homeomorphic to a compact convex subset of the Fr\'echet space 
$C_D^\infty(\bR)$, it follows again from the  Schauder--Tychonoff  Fixed Point Theorem
that $\Diff_c^\infty(\bR)$ is perfect. We refer the reader to~\cite{Epstein1984CMH} for more details.
\end{rem}

\subsection{Steps of the proof}\label{ss:mather}
Let $k,\alpha,D,E$ and $A$ be as in Theorem~\ref{t:mather-inv}.
For a function $g$, we let $\dom g$ denote the domain of $g$.
We write $\mup{\cdot}=\abss{\cdot}_{k,\alpha}$.
For a compactly supported real line homeomorphism $f$, we recall from Notation~\ref{not:If} that $I_f$ coincides with the smallest compact interval containing $\supp f$.
Let us fix the notation
\[
J:=[-2A,2A].\]

For the remainder of this section, we will construct the following.
\be[(i)]
\item\label{c:const}
Positive real numbers $C_0=C_0(k,\alpha)$ and $\delta_0=\delta_0(k,\alpha,A)$;
\item\label{c:Gamma} A map $\Gamma\co V_c^0(1)\to\Homeo_\bZ(\bR)$, which restricts to a $C^k$--continuous map \[\Gamma\restrict\co V_{J}^{k,\alpha}(\delta_0) \to \Diff_\bZ^{k,\alpha}(\bR)\]
such that for each $u\in V_{J}^{k,\alpha}(\delta_0)$ we have
\[\mup{\Gamma u} \le C_0 (|I_u| +1)\cdot \mup{u}; \]
\item\label{c:Omega} A  $C^k$--continuous map \[\Omega\co V_\bZ^{k,\alpha}(\delta_0)\to\Diff_D^{k,\alpha}(\bR)\]
such that for all $g\in  V_\bZ^{k,\alpha}(\delta_0)$ we have
 \[\mup{\Omega g} \le C_0/|D|\cdot  \mup{g},\]
and such that for all $w\in\Omega^{-1}(\dom\Gamma)$ we have
\[\Gamma\circ \Omega(w)=T(-w(0)) w;\]
\item\label{c:lambda} A diffeomorphism $\tau \in \Diff_{[-2A-1,2A+1]}^\infty(\bR)$ and a $C^k$--continuous map 
\[
\lambda\co  \left\{(u,v)\in 
V_{J}^{k,\alpha}(\delta_0)\times V_{J}^{k,\alpha}(\delta_0)
\middle| \Gamma v (\Gamma u)^{-1}\in T(\bR)\right\}
\to \Diff_{[-2A,2A+1]}^{k,\alpha}(\bR)\]
such that for all $(u,v)\in\dom\lambda$ we have
\[\tau v = \lambda (u,v) \tau u\lambda(u,v)^{-1}.\]\ee

\bp[Proof of Theorem~\ref{t:mather-inv}, assuming the above constructions]
Pick a sufficiently small $\delta>0$ so that 
\[\Psi := \Omega\circ\Gamma\co V_E^{k,\Omega}(\delta)\to \Diff_D^{k,\alpha}(\bR) \]
is a well-defined $C^k$--continuous map, and has the image inside $\dom\Gamma$.
If $g\in\dom\Psi$, then
\[
\mu(\Psi(g))=\mu(\Omega\Gamma g)\le C_0/|D|\mu(\Gamma g)\le 2C_0^2|E|/|D|\mu(g).\]
We also have that
\[
\Gamma\Psi(g) (\Gamma g)^{-1}=\left(\Gamma\Omega(\Gamma g)\right) (\Gamma g) ^{-1}
=T(-\Gamma g(0)).
\]
It follows that for 
\[
\tau \Psi(g)=\lambda \tau g \lambda^{-1}\]
for some $\lambda\in\Diff_c^{k,\alpha}(\bR)$.
This verifies the two required conditions in Theorem~\ref{t:mather-inv}.
\ep
\begin{rem}
The $C^k$--continuity of $\lambda$ was not used in the above proof.
\end{rem}

\subsection{The rolling-up process $\Gamma$}
We first define the \emph{rolling-up process} $\Gamma$, which ``mixes'' a given compactly supported homeomorphism $g$ with the translation $T$.
\begin{lem}\label{l:gam}
There uniquely exists a map 
\[\Gamma\co  
V_c^0(1)\to\Homeo_\bZ(\bR)\] such that 
for each $g\in  \dom \Gamma$, $x\in\bR$ and $r,s\in\bZ$ satisfying
\[\supp g\sse[T^{-r}(x),(Tg)^s T^{-r}(x)],\]
we have that
\[(\Gamma g)(x)=T^{r-s}(Tg)^s T^{-r}(x).\]
Furthermore, if $g$ is a $C^{k[,\alpha]}$--diffeomorphism, then so is $\Gamma g$.
\end{lem}

\bp
Let $g\in V_c^0(1)$, and put $a:=\|g-\Id\|<1$. 
Since \[\inf_{y\in\bR} |Tg(y)-y|\ge 1-a,\]
for each $x\in\bR$ we can pick $r,s\in\bZ$ such that \[\supp g\sse[T^{-r}x,(Tg)^s T^{-r}x].\]
Indeed, we can even require that \[\inf\supp g-1<x-r\le \inf\supp g\] and that \[s= \ceil*{(|I_g|+1)/(1-a)}.\]
The claim below shows that $\Gamma g\co\bR\to\bR$ is well-defined, which is to say independently of $r$ and $s$. From the claim, it follows that $\Gamma$ is uniquely determined by the given conditions.

\begin{claim*}
Let $x\in\bR$ and $r,s,r',s'\in\bZ$ satisfy
\[\supp g\sse [T^{-r}x,(Tg)^sT^{-r}x]\cap  [T^{-r'}x,(Tg)^{s'}T^{-r'}x].\]
We have that
\[T^{r-s}(Tg)^sT^{-r}(x)=T^{r'-s'}(Tg)^{s'}T^{-r'}(x).\]
\end{claim*}
By symmetry, we may assume $r-s\le r'-s'$. 
Considering both of the cases $r\ge r'$ and $r'\le r$, we easily see that
\[
(Tg)^rT^{-r}x = (Tg)^{r} T^{-(r-r')} T^{-r'}x=(Tg)^{r}(Tg)^{-(r-r')}T^{-r'}x=(Tg)^{r'}T^{-r'}x.\]
The claim follows from the computation
\begin{align*}
T^{r-s}(Tg)^sT^{-r}x&=T^{r-s}(Tg)^{s-r+r'}T^{-r'}x=T^{r-s}(Tg)^{-r+r'+s-s'}(Tg)^{s'}T^{-r'}x\\
&=T^{r-s}T^{-r+r'+s-s'}(Tg)^{s'}T^{-r'}x=T^{r'-s'}(Tg)^{s'}T^{-r'}x.
\end{align*}

We observe  $\Gamma g\in\Homeo_\bZ(\bR)$ from the computation (assuming $r\gg0$)
\[\Gamma gTx = T^{r-s}(Tg)^sT^{-r+1}x=T T^{r-1-s}(Tg)^sT^{-r+1}x=T\Gamma gx.\]
Note that positive integers $r$ and $s$ can be chosen to be constant for the definition of $\Gamma g(y)$ for all $y$ sufficiently near from a given number $x\in\bR$.
Hence, $\Gamma g$ inherits the regularity of $g$.
\ep

\begin{lem}\label{l:gam-prop1}
Let $g\in V_c^0(1)$ and let $b\in\bR$.
\be
\item Let $x\in\bR$ and $r,s\in\bZ$ satisfy $\supp g\sse [(Tg)^{-s}T^{r}(x),T^{r}(x)]$. 
Then,
\[(\Gamma g)^{-1} (x)= T^{s-r}(Tg)^{-s}T^{r}(x).\]
\item
We have that
\[
\Gamma (T(b)g T(-b))=
T(b)(\Gamma  g) T(-b).\]
\item\label{p:gam-c0}
If we set $a:=\|g-\Id\|$ and 
 \[s:=\ceil{(|I_g|+1)/(1-a)},\] then we have
\[
\|\Gamma g- \Id\|
\le sa.\]
\ee\end{lem}

\bp 
Parts (1) and (2) are immediate. Regarding part (\ref{p:gam-c0}), we note
\[ |(Tg)^s(y) -T^sy|\le 
|g (Tg)^{s-1}(y)-T^{s-1}y|
\le a+ |(Tg)^{s-1}(y)-T^{s-1}y|
\le \cdots \le sa\] for all $y\in\bR$.
Hence, for a suitable choice of $r$ we have
\[
|\Gamma gy - y|
= |(Tg)^s(y-r)+r-s - y|
= |(Tg)^s(y-r)-(y-r+s)|\le sa.\qedhere\]\ep

\begin{lem}\label{l:gam-prop2}
Let $B\in\bN$, and let $J_1,\ldots,J_B$ be unit intervals in $\bR$ satisfying $\sup J_i\le\inf J_{i+1}$ for  each $i<B$. If $h_i\in \Homeo_\bZ(\bR)$ satisfies $h_i(J_i)=J_i$ for each $i$, then \[\Gamma \left( \left(h_B\restrict_{J_B}\right)\circ  \cdots\circ \left(h_1\restrict_{J_1}\right) \right)=h_B\circ\cdots\circ h_1.\]
\end{lem}

\bp
Let us set
 \[g =\prod_B^1 h_i\restrict_{J_i},\] and
 \[h=\prod_{i=B}^1 h_i.\] 
Since $\|g-\Id\|=\max_i \|h_i\restrict_{J_i}-\Id\|<1$ we see that $\Gamma g$ is well-defined.
Without loss of generality, we may assume $J_1=[1,2]$.
Let $x\in[0,1]$. Inductively, we can find $r_i\in\bN$ for each $i<B$ such that
\[(Tg)^{r_i+r_{i-1}+\cdots+r_1+1}x= T^{r_i} h_i\cdots T^{r_2}h_2T^{r_1}h_1 Tx\in J_{i+1}.\]
Let $s=2+\sum_{i<B} r_i$. Since $h_iT=Th_i$, we deduce that
\[
\Gamma g(x) = T^{-s}(Tg)^s(x)= T^{-s}(Tg)^{s-1}(Tx)=T^{-s} T h_B \prod_{B-1}^1 T^{r_i}h_i T(x)=\prod_B^1 h_i(x).\qedhere\]
\ep

Here, the phrase ``$\delta_0=\delta_0(k,\alpha,A)$ is sufficiently small''
means that $\delta_0$ is chosen to be small, depending on the choices of $k,\alpha$ and $A$.

\begin{lem}\label{l:gam2}
If $\delta_0=\delta_0(k,\alpha,A)$ is sufficiently small, then the following hold.
\be
\item\label{p:gam-cont} 
$\Gamma$ restricts to a  $C^k$--continuous map on $V_J^k(\delta_0)$;
\item\label{p:gam-est}
Each $f\in V_J^{k,\alpha}(\delta_0)$ satisfies
\[
\mup{\Gamma f} \le 4(|I_f|+1)\cdot \mup{f}.\]
\ee\end{lem}

\bp
(\ref{p:gam-cont})
Let $f, g\in V_J^k(\delta_0)$.
The Domination Lemma implies that if $\delta_0$ is chosen to be sufficiently small, then 
$f,g\in  V_J^0(a)$ for some $a<1/100$.
Similarly as in the proof of Lemma~\ref{l:gam-prop1}
we let
\[s=\ceil*{(|J|+2)/(1-a)}.\]
Then for each $y\in[\inf J-2,\inf J]$,
we have  $\Gamma f(y)=T^{-s} (Tf)^s(y)$ and 
 $\Gamma g(y)=T^{-s} (Tg)^s(y)$.
So, \[\abss{\Gamma f-\Gamma g}_{i} \le\abss{(Tf)^{s}-(Tg)^{s}}_{i}\]
 for $0\le i\le k$.
As $\Diff_{\mathrm{ep}}^k(\bR)$ is a topological group with respect to the metric
\[ d_k(f,g)=\sup_{0\le i\le k} \abss{f-g}_i,\]
we see that $d_i((Tf)^s,(Tg)^s)$ can be made arbitrarily small by making $d_i(f,g)$ small enough. 
This proves the $C^k$--continuity of $\Gamma$.


(\ref{p:gam-est})
We may modify the definition of $s$ above as 
\[s:=\ceil*{(|I_f|+2)/(1-a)}\le 1.1(|I_f|+2)+1\le 1.1(4A+2)+1,\]
so that
\[\mup{\Gamma f}\le\mup{T^{-s}(Tf)^s}.\]

By shrinking $\delta_0$ a priori, we can find a constant $C=C(k,\alpha,A)$ as in Proposition~\ref{prop:muk-compose} such that 
\[
\mup{(Tf)^j} \le \mup{f}+(1+C\delta_0) \mup{(Tf)^{j-1}}
\]
for all $1\le j\le s$.
Shrinking $\delta_0$ further, we see that 
\[
\mup{(Tf)^s}\le \mup{f}\sum_{0}^{s-1} (1+C\delta_0)^i \le
1.1\mup{f}s
\le 1.1\mup{f}\cdot (1.1 (|I_f|+2)+1)
\le 4(|I_f|+1)\cdot\mup{f}.\qedhere\]
\ep
We have now proved part (i) in Subsection~\ref{ss:mather}.

\subsection{Detecting conjugacy from $\Gamma$}
We recall the following standard result on global existences of smooth flows.

\begin{cor}[see~\cite{LeeGTM}]\label{cor:completevf}
Every compactly supported smooth vector field on a smooth manifold without boundary uniquely admits a smooth complete flow.\end{cor}

Recall we have set $J:=[-2A,2A]$ for some $A\ge1$.
Choose an even, smooth vector field 
\[\rho\co \bR\to [0,1]\] such that $\supp\rho= [-2A-1,2A+1]$ and $\rho(J)=\{1\}$. 
By Corollary~\ref{cor:completevf}, there exists a corresponding flow 
\[ \Phi\co\bR\times\bR\to \bR\] for $\rho$. 
From this flow, we have the $0$--trajectory map and the time-$t$ map, which we record respectively as:
\[ \phi(x):=\Phi(x,0)\in \Diff_+^\infty(\bR,(-2A-1,2A+1)),\quad \tau_t(x)=\Phi(t,x)\in\Diff_{[-2A-1,2A+1]}^\infty(\bR).\]
Here, $\Diff^\infty_+(A,B)$ denote the group of orientation--preserving smooth diffeomorphisms from $A$ to $B$.
We put $\tau=\tau_1$.
It is routine to check the following.
\begin{lem}\label{l:phi}
The following conclusions hold.
\be
\item\label{p:tauphi2}
For each $ b\in\bR$, we have $\phi T_b\phi^{-1}=\tau_b\restriction_{(-2A-1,2A+1)}$.
\item\label{p:tauphi1}
$\phi\restrict_J=\Id_J$ and $|\phi(x)|\le |x|$ for $x\in \bR$.
\item\label{p:tauphi3}
Each $u\in \Homeo_J(\bR)$ satisfies  $\phi^{-1}u\phi=u$, and moreover,
 $\phi u \phi^{-1}=u\restrict_{(-2A-1,2A+1)}$.
\ee
\end{lem}

The purpose of this subsection is to prove a certain criterion for conjugacy of diffeomorphisms, which
is a crucial step in Mather's proof. 
Recall we have fixed the notations $k,\alpha,J$ and $\tau$.

\begin{lem}\label{l:rot-conj}
Let $\delta_0=\delta_0(k,\alpha,A)$ be sufficiently small.
Then there exists a $C^k$--continuous map
\[
\lambda\co\left \{
(u,v)\in V_J^{k[,\alpha]}(\delta_0)\times V_J^{k[,\alpha]}(\delta_0)\mid\Gamma v (\Gamma u)^{-1}\in T(\bR)\right\}
\to\Diff_{[-2A,2A+1]}^{k[,\alpha]}(\bR)\]
such that 
 \[ \tau v=\lambda  \tau u \lambda^{-1}.\]
\end{lem}

Thus, the diffeomorphism $\lambda$ is associated to a pair $(u,v)$ of diffeomorphisms whose images under $\Gamma$ differ by
a translation and realizes a conjugacy between $\tau v$  and $\tau u$.
We recall the notation of a sufficiently small universal constant $\epsilon_0>0$, as at the beginning of the section.
\begin{lem}\label{lem:Lamuv}
Let $\delta_0=\delta_0(k,\alpha,A)$ be sufficiently small. Then 
for each $u,v\in V_J^{k,\alpha}(\delta_0)$, there exists 
$\Lambda=\Lambda(u,v)\in\Diff_{\mathrm{ep}}^{k[,\alpha]}(\bR)$
such that  
\[\Lambda(x)=\lim_{s\to\infty}(Tv)^{s}(Tu)^{-s}(x)\]
for all $x\in\bR$.
Setting $\lambda_0:=\phi\Lambda \phi^{-1}\in\Diff_+^{k[,\alpha]}(-2A-1,2A+1)$,
we further have that 
\begin{itemize}
\item $\Lambda(x)=x$ for all $x\le -2A$;
\item $\Lambda(x)=(\Gamma v)(\Gamma u)^{-1}(x)$ for $x\ge 2A+1/2$.
\item
$
\tau v(x)
=\lambda_0\tau u\lambda_0^{-1}(x)$
for all $x\in(-2A-1,2A+1)$.
\end{itemize}
\end{lem}

Lemma~\ref{lem:Lamuv} trivially implies  $[u]=[v]$ in $H_1\left(\Diff_+^{k[,\alpha]}(\bR)\right)$. Namely, We see $[Tu]=[Tv]$ from
\[
Tv\Lambda=\lim_{s\to\infty}(Tv)^{s+1}(Tu)^{-s}
=\lim_{s\to\infty}(Tv)^{s+1}(Tu)^{-(s+1)}Tu=\Lambda Tu.\]
Since we can only draw relevant from conclusions from conjugacy in the groups of \emph{compactly supported} diffeomorphisms,
we will have to modify the support of $\Lambda$ using the flow $\Phi$.
We also emphasize that the $C^{k,\alpha}$ regularity of $\lambda_0$ at $\pm(2A+1)$ is not guaranteed yet.

\bp[Proof of Lemma~\ref{lem:Lamuv}]
If $s\in\bN$ is large enough so that 
$(Tu)^{-s}(x)\le-2A$, then 
\[(Tv)^{s+1}(Tu)^{-s-1}x=(Tv)^s (T)(T)^{-1}(Tu)^{-s}x=(Tv)^s(Tu)^{-s}x.\]
It follows that the limit exists and is realized for some $s<\infty$. 
We have  
  \[\Lambda^{-1}(x)=\lim_{s\to\infty}(Tu)^{s}(Tv)^{-s}(x).\]
Since $\Lambda$ is invertible and $s$ can be chosen to be locally constant, we see  that $\Lambda\in\Diff_+^{k[,\alpha]}(\bR)$.
By setting $s=0$, we see that $\Lambda(x)=x$ for $x\le -2A$.

Requiring  $\delta_0$ to be sufficiently small, we have that
 \[a:=d_0(\{u,v\},\Id)<1.\]
Let us pick an arbitrary $x\in [2A+1/2, 2A+3/2]$.
Setting
\[s=\ceil*{\frac{4A+3/2}{1-a}}\]
we have that \[(Tu)^{-s}(x)\le  x-(1-a)s\le 2A+3/2-(1-a)s\le -2A.\]
Shrinking $\delta_0$ further, we have that $a$ is very small compared to $A$, and that 
\[\Lambda(x) = (Tv)^s (Tu)^{-s}(x)
\ge x-(1+a)s+(1-a)s= x-2sa\ge
 2A+\frac12-2a\left(\frac{4A+3/2}{1-a}+1\right)\ge2A.\]
By the definition of $\Gamma$, we see that
\[\Gamma v\circ (\Gamma u)^{-1} (x)=  (Tv)^{s}T^{-s} \circ T^s(Tu)^{-s}(x)=\Lambda (x).\]
Moreover, for $r\in\bN$ we have that \[(Tu)^{-(s+r)}(x+r)=(Tu)^{-s}(x)\le -2A,\] and that
\[
\Lambda T^r x=(Tv)^{s+r}(Tu)^{-(s+r)}(x+r)=(Tv)^r\Lambda(x)=T^r\Lambda(x).\]
We have thus shown that \[\Lambda[2A+1/2,\infty)\sse[2A,\infty),\]
and that $\Lambda T(x)=T\Lambda(x)$ for $x\ge 2A+1/2$;
In particular, we see that $\Lambda$ is eventually periodic.

We have that  $Tv=\Lambda Tu \Lambda^{-1}$
in $\Homeo_{\mathrm{ep}}(\bR)$.
Conjugating by $\phi$, we see that
\[
\tau v=
\tau \phi v \phi^{-1}=\phi Tv\phi^{-1}=
\phi\Lambda Tu \Lambda^{-1}\phi^{-1}=
\lambda_0(\phi Tu\phi^{-1})\lambda_0^{-1}
=\lambda_0\tau u\lambda_0^{-1}.\qedhere
\]
 \ep

\bp[Proof of Lemma~\ref{l:rot-conj}]
We will follow the notations in Lemma~\ref{lem:Lamuv}.
Let $L$ be the domain of $\lambda$ as given by the hypotheses of the lemma, and let $(u,v)\in L$ so that 
\[
\Gamma v(\Gamma u)^{-1}=T_b\]
for some $b\in\bR$.
For $x\in (-2A-1,2A+1)$, we have seen $\lambda_0\tau u\lambda_0^{-1}(x)=\tau v(x)$.

\begin{claim*}
The map 
\[
\lambda_0\in \Diff_+^{k[,\alpha]}(-2A-1,2A+1)\]
 extends to a diffeomorphism 
 \[\lambda\in \Diff_{[-2A,2A+1]}^{k[,\alpha]}(\bR).\]
satisfying the following three conditions:
\[
\lambda(x)=\begin{cases}
\lambda_0=\phi\Lambda\phi^{-1}(x),&\text{ if }x\in (-2A-1,2A+1);\\
x,&\text{ if }x< -2A;\\
\Phi(b,x)=\tau_b(x),&\text{ if }x>  2A+1/2.
\end{cases}
\]
 \end{claim*}
 
Note that each of the above three functions is a local $C^{k,\alpha}$ diffeomorphism.
Hence, it remains to show that the values of $\lambda$ coincide on the overlapping intervals,
and $\supp\lambda\sse[-2A,2A+1]$.

If $-2A-1< x< -2A$, then $\phi^{-1}(x)\le -2A$ and hence,
$\phi\Lambda\phi^{-1}(x)=\phi\phi^{-1}(x)=x$. Thus, there is no conflict between the first and the second rows defining $\lambda$.
Suppose $ 2A+1/2< x<  2A+1$.
We have that $ \phi^{-1}(x)\ge x>  2A+1/2$. 
From Lemma~\ref{lem:Lamuv},
it follows that \[\phi\Lambda\phi^{-1}(x)=\phi T_b \phi^{-1}(x)=\Phi(b,x).\]
This shows that the first and the third functions coincide on the overlap.
In other words, $\lambda$ is well-defined by the above conditions. 
We note that if $x\ge  2A+1$ then
\[\Phi(b,x)=x.\]
This show that $\supp\lambda\sse[-2A,2A+1]$.

We claim that if \[(u,v),(u_1,v_1)\in \dom\lambda,\]
and if
$d_k(u,u_1)$ and $d_k(v,v_1)$ are sufficiently small, then the corresponding $\lambda(u,v)$ and $\lambda(u_1, v_1)$ will be arbitrarily $C^k$--close.
We first write $\bR$ by the union of three open intervals $I_1 , I_2, I_3$, whose closures are contained in the intervals
\[\{(-\infty, -2A),(-2A-1, 2A+1),(2A+1/2,\infty)\},\] respectively.
On $I_1$, we have \[\lambda(u,v)(x) = \lambda(u_1, v_1)(x) = x.\] On $I_2$, we have $\lambda(x) = \phi\Lambda\phi^{-1}(x)$.
Note that $\lambda(u,v)$ and $\lambda(u_1 , v_1)$ are arbitrarily $C^k$-close since $\Lambda(u,v)$ and $\Lambda(u_1 , v_1)$ are.
Specifically, by following the proof of Lemma~\ref{lem:Lamuv}, we may find an integer $s$ so that $\Lambda(u,v) = (Tv)^s (Tu)^{-s}$
and so that $\Lambda(u_1,v_1) = (Tv_1)^s (Tu_1)^{-s}$ on $[-2A-1, 2A+1]$.

For the same choice of integer $s$, we have \[b = -(2A+1) + (Tv)^s (Tu)^{-s}(2A+1),\] so that $b(u,v)$ and $b(u_1, v_1)$
are arbitrarily close. Also, the function $\Phi_0(t,x) = \Phi(t,x)-x$ is smooth and compactly supported on $I\times \bR$ for any compact interval $I$. This implies that 
\[
\lambda(u,v) - \lambda(u_1 , v_1) = \Phi_0(b(u,v), \cdot) - \Phi_0(b(u_1, v_1), \cdot)
\]
is arbitrarily $C^k$-close to 0 on $I_3$.
It follows that the map
\[(u,v)\mapsto \lambda(u,v)\]
is $C^k$--continuous. 
\ep

\subsection{The spreading process $\Omega$}
We define a process $\Omega_B$, which ``spreads'' a circle diffeomorphism to a real line diffeomorphism supported on $[-2B,2B]$. 
As described in part~(\ref{p:b2}) below, the process of spreading and rolling-up a given compactly supported diffeomorphism
returns the original diffeomorphism, up to a tranlation.

\begin{lem}\label{l:omb}
For each $B\in\bN$, there exist positive numbers $\epsilon=\epsilon(B),C=C(k,\alpha)$ and a map
\[
\Omega_B\co V_\bZ^1(\epsilon)\to \Diff_{[-2B,2B]}^1(\bR)\]
such that the following hold.
\be
\item\label{p:b1}
The map $\Omega_B$ restricts to a $C^k$--continuous map
\[
\Omega_B\restrict\co V_\bZ^1(\epsilon)\cap\Diff_+^{k[,\alpha]}(\bR)\to\Diff_{[-2B,2B]}^{k[,\alpha]}(\bR);\]
\item\label{p:b2}
 $\Gamma (\Omega_Bg)=T(-g(0))g$ for all $g\in V_\bZ^1(\epsilon)$.
\item\label{p:b3}
For some sufficiently small $\delta_0=\delta_0(k,\alpha,B)$,
each $g\in V^{k,\alpha}_\bZ(\delta_0)$ satisfies
\[\mup{\Omega_B g}\le \frac{C}{B}\mup{g}.\]
\ee
\end{lem}

The goal of this subsection is to prove the above lemma for the case $B=1$. The case $B>1$ is proved in the next subsection.
We first fix a smooth bump function $\zeta\co\bR\to[0,1]$ such that
\[\zeta(x)=\begin{cases} 1,&\text{ if }x\in\frac1{10}[-1,1];\\
0,&\text{ if }x\in\frac12+\frac1{10}[-1,1],\end{cases}\]
and such that $\zeta(x+1)=\zeta(x)+1$.
The universal constant $\epsilon_0$ will be shrunk if necessary to satisfy 
\begin{equation*}100 \epsilon_0\le 1/ (1+\abss{\zeta}+\abss{\zeta'}).\end{equation*}

\begin{lem}\label{l:disj-sup}
Let $J_1,\ldots,J_B$ be compact intervals in $\bR$ with disjoint interiors.
If $g_i\in\Diff_{J_i}^{k,\alpha}(\bR)$ for each $i$,
then
\[
\mup{g_B\cdots g_1} \le 2\max_j \mup{g_j} .\]
\end{lem}

\bp
We set $g =g_B\cdots g_1$ and $N= \max_i \mup{g_i} $.
For $x>y$, we will estimate $\der{g}{k}x-\der{g}{k}y$.
If $x,y\not\in\cup_i J_i$, then $\der{g}{k}x=\der{g}{k}y=\der{\Id}{k}$. So, we assume $x\in J_i$ for some $i$.

If $y\in J_i$, then \[\abs{\der{g}{k}x-\der{g}{k}y}=\abs{\der{g_i}{k}x-\der{g_i}{k}y}\le N\alpha(x-y).\]
If $y\in J_j$ for some $j\ne i$, then we can find $x'\in \partial J_i$ and $y'\in\partial J_j$ such that 
\begin{align*}
|x-y|&= |x-x'|+|x'-y'|+|y'-y|,\\
|\der{g}{k}x-\der{g}{k}y|&\le |\der{g}{k}x - \der{g}{k}x'|+|\der{g}{k}y'-\der{g}{k}y|\le N(\alpha(x-x')+\alpha(y-y'))\le 2N\alpha(x-y).\end{align*}
The case when $y\not\in\cup_j J_j$ is similar.\ep

It will be convenient for us to employ the following concise notation.
\begin{notation}
Let $f,g\co U\to\bR$ be maps on some set $U$ (possibly, a function space). We say 
\[
f\le O(g)\]
if there exist positive numbers $C$ and $\epsilon$ such that $f(u)\le Cg(u)$ 
for all $u\in U$ satisfying $|g(u)|\le\epsilon$.
\end{notation}
Note that $f\le O(|g|)$ is equivalent to the combination of the following two conditions:
\begin{itemize}
\item whenever $g(u)=0$, we have $f(u)\le 0$;
\item we have \[\lim_{\epsilon\searrow0}\sup\left\{\frac{f(u)}{|g(u)|}\co u\in g^{-1}[-\epsilon,\epsilon]\text{ and }g(u)\ne0\right\}<\infty.\]
\end{itemize}

\bp[Proof of Lemma~\ref{l:omb} for the case $B=1$]
Let $g\in V^1_\bZ(\epsilon_0)$ be given. We set
  \[h := T(-g(0))g,\quad
  h_0 := \zeta h+(1-\zeta)\Id,\quad h_1 := h\circ (h_0)^{-1}.\]
Then we let
\[\Omega_1g:= \left(h_1\restrict_{[1,2]}\right)\circ  \left(h_0\restrict_{[-3/2,-1/2]}\right)
\in\Diff_+^1(\bR).\]

Since $h[0,1]=[0,1]$, we see that $\abss{h-\Id}\le \abss{h-\Id}_1\le\epsilon_0$. From
\begin{align*}
h_0(x+1)&=\zeta(x) (h(x)+1) + (1-\zeta(x))(x+1)=h_0(x)+1,\\
\abss{h_0-\Id}_1&\le\abss{\zeta' (h-x)+\zeta (h'-1)}\le  \left(\abss{\zeta'}+\abss{\zeta}\right)
\abss{h-\Id}_1\le 1/10,
\end{align*}
we have that $h_0,h_1\in\Diff_\bZ^1(\bR)$. We estimate (with generous bounds) that
\begin{align*}
\abss{h_1-\Id}_1&\le\abss{(h'/ h_0')\circ h_0^{-1}-1}\le \abss{h'-h_0'}/\inf h_0'\\
&\le 2\abss{h-h_0}_1\le 2(1+\abss{\zeta}+\abss{\zeta'})\abss{h-\Id}_1\le1/10.
\end{align*}
Succinctly, we may write that 
\[\abss{h_0'-1},\abss{h_1'-1}\le O(\abss{h'-1}).\]
Here,  we are considering the real valued maps defined on $V^1_\bZ(\epsilon_0)$:
\[
g\mapsto \abss{h'-1},\quad g\mapsto\abss{h_0'-1},\quad g\mapsto\abss{h_1'-1}.\]

Note that 
\[h_0\left(\frac1{10}[-1,1]\right)=h\left(\frac1{10}[-1,1]\right)\supseteq  \frac{[-1,1]}{20}.\]
It follows that  $h_1(x)=x$ for \[x\in \frac{[-1,1]}{20}.\] It is obvious that $h_0(x)=x$ for \[x\in\frac12+\frac{ [-1,1]}{20}.\] 
Puting \[J_0=[-3/2,-1/2],\, J_1=[1,2]\] and setting $g_i=h_i\restrict_{J_i}$ for $i=0,1$.
We then have that  $\Omega_1(g)=g_1g_0\in \Diff_{[-2,2]}^1(\bR)$.

For part (\ref{p:b1}), we first note that the regularity of $h$ inherits that of $g$. 
Since the correspondence $g\mapsto (h,h_0,h_1)$ is $C^k$--continuous, so is the map $\Omega_1$.

Part  (\ref{p:b2}) follows from
Lemma~\ref{l:gam-prop2}, together with the computation
\[
\Gamma\Omega_1(g)=\Gamma(g_1g_0)=h_1 h_0 = h=T(-g(0))g.\]
For part (\ref{p:b3}), 
we require $V_\bZ^{k,\alpha}(\delta_0)\sse V_\bZ^1(\epsilon_0)$.
Using estimates of $\mup{\cdot}$ under taking multiplications and inverses,
and shrinking $\delta_0$ if necessary we have
\begin{align*}
\mup{h_0}&
\le
C\sum_{0\le i\le k}\brac*{\der{\zeta}{k-i}\der{(h-\Id)}{i}}_\alpha
\le
C'\sum_{0\le i\le k}\left(\abss{h-\Id}_{i,\alpha}+\abss{h-\Id}_i\right)
\le O\left(\mup{h}\right),\\
\mup{h_1}&\le \mup{h}+\mup{h_0^{-1}}
+C\cdot\mup{h}\cdot\mup{h_0^{-1}}
\le O\left(\mup{h}\right).
\end{align*} 
for some constant $C,C'$ depending on $k,\alpha$ and $\zeta$.
Here, the notation $O(\cdot)$ is used in the context of real-valued maps on $V_\bZ^{k,\alpha}(\delta_0)$. The Domination Lemma can be applied for $i=0$ since $h(0) = 0$.

Note that $\mup{g_i}\le\mup{h_i}$ for $i=0,1$. 
By Lemma~\ref{l:disj-sup}, we have
\[\mup{\Omega_1(g)}=
\mup{g_1 g_0}
\le2\max(\mup{g_1 },\mup{g_0})\le
O(\mup{h})=O(\mup{g}).
\qedhere\]
\ep

\subsection{Discrete isotopy}
The case $B>1$ of Lemma~\ref{l:omb} will require a ``discretized form'' of an isotopy, 
as described here. 
We first note the subadditivity of $[\cdot]_\alpha$, made precise as follows.
\begin{lem}\label{l:conv}
For a pointwise convergent series of real functions $\sum_i f_i$,
we have \[\brac*{\sum_i f_i}_\alpha\le\sum_i \brac*{f_i}_\alpha.\]
\end{lem}
\bp
We may certainly assume the right  hand side is finite. 
For $F_N=\sum_{1}^N f_i$, we have
\[
|F_Nx-F_Ny|\le \sum_1^N |f_ix-f_iy|\le\sum_1^N [f_i]_\alpha\alpha(x-y).\]
Put $f=\sum_i f_i$. By sending $N\to\infty$, we obtain the desired inequality
\[|fx-fy|\le \sum_1^\infty[f_i]_\alpha \alpha(x-y).\qedhere\]
\ep

We will also need the following technical estimate.
\begin{lem}\label{l:kal}
Fix $\lambda\in[0,1]$,
and write \[\abs{\cdot}=\abss{\cdot-\Id}_{k[,\alpha]}.\]
For $u\in \Diff_\bZ^{k[,\alpha]}(\bR)$ satisfying $u0=0$
and
sufficiently $C^{k[,\alpha]}$--close to $\Id$,
 we have
\[\abs*{u\circ (\lambda u+(1-\lambda)\Id)^{-1}}\le (1-\lambda)\abs{u}+O(\abs{u}^2).\]
\end{lem}
We will set $v = \lambda u+(1-\lambda)\Id$. The conclusion of Lemma~\ref{l:kal} is simply that
\[
\lim_{\epsilon\searrow0}
\sup\left\{
\frac{\abs{u\circ v^{-1}}-(1-\lambda)\abs{u}}{\abs{u}^2}
\co
u\in V_\bZ^{k[,\alpha]}(\epsilon)\text{ s.t. }u0=0
\right\}<\infty.
\]

\bp[Proof of Lemma~\ref{l:kal} in the case $\abs{\cdot}=\abss{\cdot-\Id}_{1}$]
We may first require $\abss{u'-1}\le \epsilon_0<1/4$. Write $V=v^{-1}$. 
Since \[v' = 1 + \lambda (u'-1)\ge 3/4,\] we see that $v\in \Diff^1_\bZ(\bR)$ and  
\[ \abss{V'-1}\le\frac{\abss{v'-1}}{\inf v'}\le 2\abss{v'-1}=2\lambda\abss{u'-1}.\]
We have a power series expansion (for $|\lambda z|<1$) 
\[(1+z)/(1+\lambda z) = 1+(1-\lambda)z +\sum_{i\ge2} a_i z^i.\]
Here,  $|a_i|=(1-\lambda)\lambda^{i-1}$.
By setting $z = u'\circ V-1$ and
\[C_1=\sup_{\lambda\in[0,1]}\sum_{i\ge2} |a_i|(1/4)^{i-2} = \sup_{\lambda\in[0,1]} \lambda(1-\lambda)\sum_{i\ge2} (\lambda/4)^{i-2}<\infty,\]  we see that
\begin{align*}
(u \circ V)' &= (u'/v')\circ  V  = 1+(1-\lambda)(u'\circ V-1)+\sum_{i\ge 2} a_i (u'\circ V-1)^i,\\
\abss*{\left(u \circ V\right)'-1} &
\le (1-\lambda)\abss{u'-1}
+C_1\abss{u'-1}^2.\qedhere
\end{align*}
\ep

\bp[Proof of Lemma~\ref{l:kal} in the case $\abs{\cdot}=\abss{\cdot-\Id}_{1,\alpha}$]
Using Lemma~\ref{l:conv} and continuing from the proof in the case $\abs{\cdot}=\abss{\cdot-\Id}_{1}$, we have
\[
\abs{u \circ V}=\brac*{(u \circ V)'}_\alpha
\le
(1-\lambda)\brac*{u'\circ V}_\alpha+\sum_{i\ge 2} |a_i|\brac*{ (u'\circ V-1)^i}_\alpha.
\]
From Lemmas~\ref{l:dom} and~\ref{l:derivation}, we have that $\abss{u'-1}\le O(|u|)$ and that
\[
\brac*{u'\circ V}_\alpha  
\le
\abs{u}(1+\abss{V-\Id}_1)
\le
\abs{u}+O(|u|^2).\]
Thus, we have that
\begin{align*}
\sum_{i\ge2} |a_i|\cdot\brac{(u'\circ V-1)^i}_\alpha
&\le \sum_{i\ge2} i|a_i| \cdot\abss{u'\circ V-1}^{i-1}\cdot\brac*{u'\circ V-1}_\alpha\\
&\le \left(\sum_{i\ge2} i|a_i| \cdot\abss{u'-1}^{i-2}\right) O( |u|^2)=O(|u|^2).\qedhere\end{align*}
\ep

We will not actually use the following two cases, but we include the proofs for completeness; Mather needed these cases in order to prove the simplicity of $\Diff_c^k(\bR^n)$ for $n\ge k> 1$.
\bp[Proof of Lemma~\ref{l:kal} in the case $\abs{\cdot}=\abss{\cdot-\Id}_{k}$ for $k>1$]
There exist  indices $\gamma=(j_1,\ldots,j_i)$ and corresponding constants $C_\gamma>0$ such that 
\begin{equation}
\der{(u\circ V)}{k}=
(\der{u}{k}\circ V )(V')^k+(u'\circ V)\der{V}{k}+\sum_\gamma C_\gamma (\der{u}{i}\circ V)\prod_{t=1}^i\der{V}{j_t}.\end{equation}
Setting $j_0=1$, we also saw Lemma~\ref{l:inverse} that
\[\der{V}{k}=-(\der{v}{k}\circ V)(V')^{k+1}-\sum C_\gamma (\der{v}{i}\circ V)\prod_{t=0}^i\der{V}{j_t}.\]
Applying the same computation for $1\le i\le k$, we see that \[\abss{\der{V}{i}}\le O(\abss{\der{v}{i}})\le O(|u|).\] We conclude
\begin{align*}
\abss{u\circ V}_k
&\le \abss*{(\der{u}{k}\circ V)(V')^k -(u'\circ V)(\der{v}{k}\circ V)(V')^{k+1}}+O(|u|^2)\\
&=\abss*{\frac{
\left(\der{u}{k}\circ V\right)\left(1+\lambda (u'\circ V-1)\right)
-
\left(u'\circ V\right)\left(\lambda \der{u}{k}\circ V\right)
}{(1+\lambda(u'\circ V-1))^{k+1}}}+O(|u|^2)\\
&=(1-\lambda)\abss*{\frac{\der{u}{k}\circ V}{(1+\lambda(u'\circ V-1))^{k+1}}}+O(|u|^2)
\le(1-\lambda)|u|+O(|u|^2)\qedhere.\end{align*}\ep

\bp[Proof of Lemma~\ref{l:kal} in the case $\abs{\cdot}=\abss{\cdot-\Id}_{k,\alpha}$ for $k>1$]
Whenever $1<i\le k$ we still have 
\[\abss{V}_{i[,\alpha]}\le O\left(\abss{v}_{i[,\alpha]}\right)\le O(|u|).\] 
We also note
\[\brac*{\der{u}{i}\circ V}_\alpha\le\brac*{\der{u}{i}}_\alpha(1+\abss{V'-1})\le |u|+O(|u|^2).\]
Using the same computation as in the case $\abs{\cdot}=\abss{\cdot-\Id}_{k}$, we have
\[
\brac*{\der{(u\circ V)}{k}}_\alpha\le
(1-\lambda)\abs*{\frac{\der{u}{k}\circ V}{(1+\lambda(u'\circ V-1))^{k+1}}}+O(|u|^2)
\le(1-\lambda)|u|+O(|u|^2).\qedhere\]\ep

\bd[Discrete Isotopy]
Let $B\in\bN$.
For each $i=0,1,\ldots,B$, 
we define the $i^{th}$ \emph{$B$--step discrete isotopy map}
\[\operatorname{DI}_B^i\co \{h\in \Diff_+^1(\bR)\co \abss{h-\Id}_1<1\} \to\Diff_+^1(\bR)\]
as follows.
Let $h\in\Diff_+^1(\bR)$ satisfy $\abss{h-\Id}_1<1$.
Consider the isotopy
\[\Phi(t,x)=t h(x)+(1-t)x,\quad t\in[0,1].\]
We set \[g_i(x)=\Phi(i/B,x),\, h_i=g_i\circ g_{i-1}^{-1},\]
and define
$\operatorname{DI}_B^i(h) =h_i$.
\ed
Note that $h=\prod_B^1\operatorname{DI}_B^i(h)$.
\begin{lem}\label{l:d-iso}
Let $B\in\bN$.
If $h\in \Diff_\bZ^{k[,\alpha]}(\bR)$ satisfies that $h0=0$ and 
is sufficiently $C^{k[,\alpha]}$--close to $\Id$, then
 \[\abss*{\operatorname{DI}^i_B(h)}_{k[,\alpha]}\le\frac1B \abss*{h}_{k[,\alpha]}+O\left(\abss*{h}_{k[,\alpha]}^2\right),\quad i=1,2,\ldots,B.\]
\end{lem}
\bp
Let $\abs{f}:=\abss{f-\Id}_{k[,\alpha]}$.
Set $g_0:=\Id$ and 
\[g_i:=\frac{i}B h + \left(1-\frac{i}{B}\right)\Id=\frac1{B-i+1}h+\left(1-\frac1{B-i+1}\right)g_{i-1}\]
for $i\ge1$. As above, we let
\[h_i:=g_i g_{i-1}^{-1}=\frac1{B-i+1}h g_{i-1}^{-1}+\left(1-\frac1{B-i+1}\right)\Id.\]

Applying  Lemma~\ref{l:kal} with $u=hg_{i-1}^{-1}$ and $\lambda = 1/(B-i+1)$,
we see
\[\abs{h g_i^{-1}}=\abs{h g_{i-1}^{-1} (g_{i-1}g_i^{-1})} \le \left(1-\frac1{B-i+1}\right)\abs{h g_{i-1}^{-1}}+O(\abs{h g_{i-1}^{-1}}^2).\]
An easy induction shows that \[\abs{h g_i^{-1}}\le (1-i/B)|h|+O(\abs{h}^2).\] It thus follows that
\[
\abs{h_i} =\frac{1}{B-i+1}\abs{h g_{i-1}^{-1}}\le \frac{1}{B}\abs{h}+O(\abs{h}^2).\qedhere\]
\ep

We can now define the spreading process $\Omega_B$ for a general $B\in\bN$.
Conceptually, if $g(0)=0$ then $\Omega_B(g)$ is obtained by applying $\Omega_1$ to the $i$--th $B$--step discrete isotopy from $\Id$ to $g$ on the $i$--th subinterval of $[-2B,2B]$.

\bd\label{d:B-frag}
For each $B\in\bN$ and for a sufficiently small $\epsilon>0$, we define
\begin{align*}
\Omega_B&\co V_\bZ^1(\epsilon)\to\Diff_{[-2B,2B]}^1(\bR),\\
\Omega_Bg&=
\prod_{i=B}^1
{T(-2B-2+4i)}\circ\Omega_1\operatorname{DI}^i_B\left(T(-g0)g\right)\circ
{T(2B+2-4i)}.\end{align*}
\ed

Finally, we complete the proof of Lemma~\ref{l:omb} as follows.
Part (1) is clear. Part (2) is a consequence of Lemma~\ref{l:gam-prop2}.
The estimate of part (\ref{p:b3}) follows from Lemmas~\ref{l:disj-sup} and~\ref{l:d-iso},
as well as our proof for the case $B=1$.

It is now immediate to see that the choices
\begin{itemize}
\item $B=1$ for $k\ge2$;
\item $B=A$ for $k=1$,
\end{itemize}
together with the map $\Psi = \Omega_B\circ\Gamma$ satisfy all the
necessary conditions in Subsection~\ref{ss:mather}. Theorem~\ref{t:mather-inv} is thus established.

\section{Critical regularity of groups of homeomorphisms}\label{sec:critical}

In this section, we wish to give and account of the state of the art of the \emph{critical regularity problem} for subgroups of $\Homeo_+(M)$, where $M\in\{I,S^1\}$ with a view towards Mather's results.
Due to the technically complicated nature of the subject, we do not pretend to give a self-contained account, and we only offer the reader a general overview of the statements of the main results. A precise formulation of the problem is as follows:

\begin{prob}\label{prob:crit-reg}
Let $G\leq\Homeo_+(M)$ be a subgroup. What is the value of \[\CR_M(G)=\sup\{r\mid G\textrm{ is a subgroup of }\Diff_+^{r}(M)\}?\]
\end{prob}

Here, we are asking that $G$ be realized as an abstract subgroup of $\Diff_+^{r}(M)$, so that a ``smoothing" of $G$ may not be topologically conjugate to the original realization of $G$ as a group of homeomorphisms. We will mostly be concerned with the case where $G$ is a countable subgroup of $\Homeo_+(M)$, and we will give some applications of the theory for countable subgroups to continuous groups as well.

Observe that since it is defined as a supremum, $\CR_M(G)$ may or may not be realized.

\subsection{Critical regularity of specific groups}
The starting point for investigating critical regularity of groups if a result of Deroin--Kleptsyn--Navas~\cite{DKN2007}, which shows that a countable subgroup of $\Homeo_+(M)$ is always topologically conjugate to a group of bi--Lipschitz homeomorphisms. Bi--Lipschitz homeomorphisms are differentiable almost everywhere, and it is reasonable to consider $1$ to be the absolute lower bound on the critical regularity of a countable group of homeomorphisms. There are, however, examples of countable groups of homeomorphisms which are not realized by diffeomorphisms~\cite{Calegari2006TAMS}, which in higher regularity leads to the rather subtle ambiguity of whether $(k)+1=(k+1)$, i.e. whether a countable group of $C^k$ diffeomorphisms whose $k^{th}$ derivatives are Lipschitz is realized as a group of $C^{k+1}$ diffeomorphisms. See Subsection~\ref{ss:open} below.

Generally speaking, Problem~\ref{prob:crit-reg} is extremely difficult and practically inaccessible with current technology for nearly all groups
with finite critical regularity.
For countable groups of diffeomorphisms, there are many groups for which the critical regularity is easily shown to be infinite (such as abelian groups, free groups, etc.) whereas finite critical regularities are only known for three essentially different examples.

For one, Navas~\cite{Navas2008GAFA}
proved that the Grigorchuk-Machi group of intermediate growth has critical regularity exactly $1$ and that this regularity is realized.

The other two groups for which the critical regularity is known are both nilpotent. On the one hand, the
Plante--Thurston Theorem~\cite{PT1976} implies that nilpotent groups have critical regularity at most $2$ (in fact they cannot be realized as groups of diffeomorphisms whose derivatives are Lipschitz). On the other hand, Farb--Franks~\cite{FF2003} and
Jorquera~\cite{Jorquera} proved that finitely generated torsion--free nilpotent groups can be realized as groups of diffeomorphisms for any $M\in\{I,S^1\}$. Castro--Jorquera--Navas~\cite{CJN2014} proved that the integral Heisenberg group $H$ has critical regularity exactly $2$.
By the Plante--Thurston Theorem, this means that $H$ can be realized as a subgroup of $\Diff^{2-\epsilon}(M)$ for every $\epsilon>0$, and the critical regularity is not realized. Jorquera--Navas--Rivas~\cite{JNR2018} proved that the group $N_4$ of $4\times4$ integral unipotent matrices has critical regularity exactly $3/2$ on the interval. It is not known whether or not this critical regularity is realized or not.

\subsection{Higher regularity and other properties of groups}

Higher regularity has several deep and mysterious relations with other properties of groups. For one, higher regularity of group actions is incompatible with Kazhdan's property (T). Navas~\cite{Navas2002ASENS},
building on work of Witte, Ghys, and Burger--Monod~\cite{Witte1994,Ghys1999,BM1999} proved that if $G$ is an infinite group with property (T) then $G$ is not contained in $\Diff^1(I)$, and that for all $\epsilon>0$ we have that $G$ is not a subgroup of $\Diff^{1.5+\epsilon}(S^1)$. Bader--Furman--Gelander--Monod~\cite{BFGM2007AM} improved this latter estimate to show that $G$ is not a subgroup of $\Diff^{1.5}(S^1)$. Thus, property (T) places \emph{a priori} restrictions on the critical regularity of a group of homeomorphisms.

Second, higher regularity of group actions puts nontrivial restrictions on the geometry of a group. Specifically,
Navas~\cite{Navas2008GAFA} proved that the if $G<\Diff^{1+\epsilon}(I)$ for $\epsilon>0$ then $G$ either has polynomial growth or exponential growth. In particular, intermediate growth of a group places \emph{a priori} restrictions on the critical regularity of a group of homeomorphisms.

\subsection{Compactness versus non--compactness and the M\"uller--Tsuboi flattening process}
Recall that we defined 
\[
\Diff_I^{k[,\alpha]}(\bR) = \{f \in \Diff_+^{k[,\alpha]}(\bR)\mid \supp f \sse I\}.\]
We clearly have that $\Diff_I^{k[,\alpha]}(\bR)\sse \Diff_+^{k[,\alpha]}(I)$.
It is less clear that embeddings exist in the other direction.

\begin{thm}\label{thm:embeddability}
There exists a $C^\infty$--homeomorphism of $I$ conjugating 
$\Diff_+^{k[,\alpha]}(I)$ into $\Diff_I^{k[,\alpha]}(\bR)$.\end{thm} 

M\"uller and Tsuboi (independently) discovered the $C^k$--cases~\cite{Tsuboi1984}, which was extended to $C^{k,\alpha}$--cases by the first two authors.
In particular, the spectra of isomorphism types of
finitely generated subgroups coincide for $\Diff_+^{k[,\alpha]}(I)$ and $\Diff_I^{k[,\alpha]}(\bR)$. We refer the reader to the appendix of~\cite{KK2017crit} for the proof of Theorem~\ref{thm:embeddability}.

\subsection{Technical tools and critical regularity ranges}

In this subsection, we give two $C^1$--obstructions for group actions which were developed by Baik--Kim--Koberda and Kim--Koberda in their investigation of the critical regularity of right-angled Artin groups and mapping class groups of surfaces.

\subsubsection{The Two-jumps Lemma}
The Two-jumps Lemma was developed by Baik--Kim--Koberda~\cite{BKK2016}, and very roughly says that two diffeomorphisms cannot pull points out of each others' supports on infinitely many small scales. The Two-jumps Lemma rules out a degree of what could be called \emph{oscillatory behavior} for pairs of diffeomorphisms. We remark that the difficulty of the result is primarily in formulating the correct statement, whereas the proof is a relatively straightforward application of the Mean Value Theorem. We state here a version which appears as Lemma 6.3 of~\cite{BKK2016}, which is easier to parse than the most general statement:

\begin{thm}[Two--jumps Lemma,~\cite{BKK2016}]\label{l:fg}
Let $M\in\{I,S^1\}$ and let $f,g\co M\to M$ be continuous maps.
Suppose  $(s_i), (t_i)$ and $(y_i)$ are infinite sequences of points in $M$
such that for each $i\ge1$, one of the following two conditions hold:
\be[(i)]
\item
$f(y_i)\le s_i = g(s_i) < y_i < t_i = f(t_i) \le g(y_i)$;
\item
$g(y_i)\le t_i = f(t_i) < y_i < s_i = g(s_i) \le f(y_i)$.
\ee
If $|g(y_i)-f(y_i)|$ converges to $0$ as $i$ goes to infinity,
then $f$ or $g$ fails to be $C^1$.
\end{thm}

\begin{figure}[h!]
  \tikzstyle {bv}=[black,draw,shape=circle,fill=black,inner sep=1pt]
\begin{tikzpicture}[>=stealth',auto,node distance=3cm, thick]
\draw  (-4,0) node (1) [bv] {} node [below]  {\small $f(y_i)$} 
-- (-2,0)  node  (2) [bv] {} node [below] {\small $s_i$}
-- (0,0) node (3) [bv] {} node [below] {\small $y_i$} 
-- (2,0)  node (4) [bv] {} node [below] {\small $t_i$}
-- (4,0) node (5) [bv] {} node [below]  {\small $g(y_i)$};
\path (3) edge [->,bend right,red] node  {} node [above] {\small $f$}  (1);
\path (3) edge [->>,bend right,blue] node  {} node [below] {\small $g$} (5);
\draw [->>, blue] (2)  edge [out = 200,in=-20,looseness=50] node [below] {\small $g$} (2);
\draw [->, red] (4)  edge [out = 20,in=160,looseness=50] node [above] {\small $f$} (4);
\end{tikzpicture}%
\caption{The Two--jumps Lemma.}
\label{f:fg}
\end{figure}
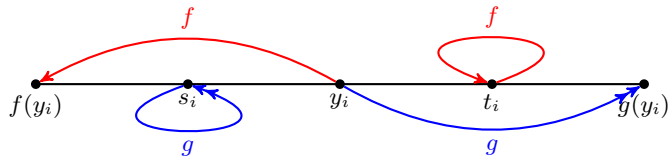

Theorem~\ref{l:fg} is illustrated in Figure~\ref{f:fg}.

The usefulness of the Two-jumps Lemma may not be immediately clear to the reader. To suggest how it might be applied, the reader may consider a group $G$ endowed with a hypothetical faithful action on $M$ of some sufficiently high regularity. Oftentimes, the algebraic structure of $G$ can be exploited to force any faithful action of $G$ to admit elements whose actions satisfy the hypotheses of the Two-jumps Lemma, whence one concludes that the hypothetical action did not exist. The groups whose algebraic structure has been most amenable to this strategy are ones which admit a complicated pattern of commuting elements and nonabelian free groups, such as is found in many right-angled Artin groups.

\subsubsection{The $\langle a,b,t\rangle$ Lemma}

The $\langle a,b,t\rangle$ Lemma was developed by Kim--Koberda~\cite{KKFreeProd2017}, and strictly speaking it is a consequence of the Two-jumps Lemma. Its usefulness will be immediately clearer to the reader, since the hypotheses are straightforward to parse and are (in principle) easy to check. As opposed to the Two-jumps Lemma, the $\langle a,b,t\rangle$ Lemma has a much more involved proof.

\begin{thm}[The $\langle a,b,t\rangle$ Lemma]\label{thm:abt}
Let $\{a,b,t\}\sse\Diff^1(M)$, and suppose that \[\supp a\cap\supp b=\varnothing.\] Then the group $\langle a,b,t\rangle$ is not isomorphic to $\bZ^2*\bZ$.
\end{thm}

\subsubsection{Applications}

The Two-jumps Lemma and the $\langle a,b,t\rangle$ Lemma cannot, as far as we can tell, be used to compute exact values of critical regularities of groups. However, they can be used to give upper bounds on the critical regularities of large classes of groups. We will state precisely the conclusions one can draw from right-angled Artin groups and mapping class groups of surfaces.

\begin{thm}[See~\cite{KKFreeProd2017}]\label{thm:cr-raag}
Let $A(\gam)$ be a right-angled Artin group. Then exactly one of the following conclusions holds:
\begin{enumerate}
\item
The group $A(\gam)$ decomposes as a (possibly trivial) direct product of groups, each of which is a (possibly trivial) free product of free abelian groups. In this case, we have $\CR_M(A(\gam))=\infty$;
\item
We have $\CR(A(\gam))\in [1,2]$.
\end{enumerate}
\end{thm}

One of the original motivations behind the development of Theorem~\ref{thm:cr-raag} was the study of mapping class group actions on one--manifolds, which has a very rich history (see~\cite{FF2001,Farb2006}). Here, $S$ will be an orientable finite-type surface with genus $g$ and a total of $n$ punctures and boundary components. By analogy to Theorem~\ref{thm:cr-raag}, we have the following result:

\begin{thm}[See~\cite{BKK2016}]\label{thm:cr-mcg}
Let $\Mod(S)$ be the mapping class group of an orientable surface $S$. Then exactly one of the following conclusions holds:
\begin{enumerate}
\item
There is a finite index subgroup $G\leq\Mod(S)$ such that $\CR_M(G)=\infty$. In this case, we have $3g-3+n\leq 1$ and the surface is sporadic;
\item
If $G\leq\Mod(S)$ is a finite index subgroup which lies in $\Homeo_+(M)$ then $\CR_M(G)\in [1,2]$.
\end{enumerate}
\end{thm}

In the second conclusion of Theorem~\ref{thm:cr-mcg}, we need to assume that the group $G$ lies in $\Homeo_+(M)$ since we do not know if mapping class groups always admit such subgroups. In particular, if $n=0$ and $g\geq 2$ then it is unknown if mapping class groups are virtually linearly (or even cyclically) orderable.

\subsection{The critical regularity spectrum}

We now consider the converse of Problem~\ref{prob:crit-reg}, which is to say we start with an $\alpha\in[1,\infty)$ and attempt to produce a countable group with that critical regularity. The set of values $\alpha$ for which this is possible is called the \emph{critical regularity spectrum} of $\Homeo_+(M)$. It turns out that one can do much more than just populate the whole critical regularity spectrum.

For a concave modulus $\alpha$ or for $\alpha\in\{0,\mathrm{bv}\}$, 
the set of all countable subgroups of $\Diff_+^{k,\alpha}(M)$ is denoted as $\GG^{k,\alpha}(M)$. The strongest succinct result as obtained by Kim--Koberda~\cite{KK2017crit} can be stated as follows: 

\begin{thm}\label{thm:spectrum}
For each $k\in\bN$,
and for each concave modulus $\mu\gg\omega_1$, there exists a finitely generated group  
$Q=Q(k,\mu)\le\Diff_+^{k,\mu}(I)$ such that the following hold.
\be[(i)]
\item
 $[Q,Q]$ is simple
 and every proper quotient of $Q$ is abelian;
\item
if $\alpha=\mathrm{bv}$, or if $\alpha$ is a concave modulus satisfying $\mu\gg \alpha\succ_k 0$, then
\[
[Q,Q]\not\in
\GG^{k,\alpha}(I)\cup \GG^{k,\alpha}(S^1).
\]
\ee
\end{thm}

Here, the notation $\alpha\succ_k 0$ means that for some $\delta>0$ we have \[\lim_{t\to +0}\sup_{0<x<\delta}t^{k-1}\alpha(tx)/\alpha(x)=0.\] 
This is equivalent to saying that either $k>1$ or $\alpha$ is sub-tame.
The notation $\alpha\ll\mu$ means \[\lim_{x\to+0}\frac{\alpha(x)\log^K(1/x)}{\mu(x)}=0\] for all $K>0$.
We obtain the following immediate consequence:

\begin{cor}\label{cor:spectrum}
For all $\alpha\in [1,\infty)$, there exists a finitely generated group and a countable simple group with critical regularity exactly $\alpha$. Moreover, for all $\alpha$ there exist examples for which the critical regularity is achieved and for which the critical regularity is not achieved.
\end{cor}

The proof of Theorem~\ref{thm:spectrum} is very complicated, and to our knowledge even the weaker statement of Corollary~\ref{cor:spectrum} does not admit a proof which is essentially simpler than that of the full strength result. To give the reader a vague sketch, we note that the $\langle a,b,t\rangle$ Lemma is one of the seeds which makes Theorem~\ref{thm:spectrum} plausible. Aside from that, one needs to develop several further tools. One is a dynamical interpretation of regularity of diffeomorphisms, combined with a recipe for producing diffeomorphisms of a given regularity. Then, one needs a combinatorial apparatus to make sense of certain groups of diffeomorphisms move points faster than any ``comparable" group of diffeomorphisms with strictly higher regularity. Finally, one analyzes the supports of elements to show that a certain tailor-made group of diffeomorphisms admits no faithful representations in strictly higher regularity.

\subsection{Applications to continuous groups}

We now connect critical regularity back to the algebraic structure of full diffeomorphism groups, and in particular Mather's results. Critical regularity gives control over the existence of nontrivial homomorphisms between diffeomorphism groups. In these situations, one often encounters the phenomenon of \emph{automatic continuity}. A result of Mann says that if $X$ is a compact manifold and if $\Homeo_0(X)$ denotes the group of homeomorphisms of $X$ which are isotopic to the identity then any homomorphism from $\Homeo_0(X)$ to a separable topological group is automatically continuous (see~\cite{Mann2016GT}). In particular, this can be used to show that $\Homeo_0(X)$ does not admit nontrivial homomorphisms into diffeomorphisms of higher regularity (i.e. the full group $\Homeo_0(X)$ admits no \emph{algebraic smoothings}).

Recall we have set
\[
\DD(M;k,\alpha):=\begin{cases}
\Diff_c^{k,\alpha}(\bR)&\text{ if }M=\bR,\\
\Diff_+^{k,\alpha}(S^1)&\text{ if }M=S^1.\end{cases}\]
In the notation of Theorem~\ref{thm:spectrum}, if $\mu\gg\alpha\succ_k 0$ then we may view elements of $\DD(M;{k,\alpha})$ as elements of $\DD(M;{k,\mu})$, since $\alpha$ a stronger modulus of continuity than $\mu$. Thus, there are obvious injective homomorphisms $\DD(M;{k,\alpha})\to\DD(M;{k,\mu})$. As for maps going the other way, there are almost no abstract homomorphisms at all.

Combining Theorem~\ref{thm:spectrum} with Mather's Theorem~\ref{thm:mather-1d}, we have the following result (cf. Corollary~\ref{cor:cont-hom}):

\begin{thm}\label{thm:homomorphism}
Let $k\in\bN$ and let $\mu\gg\alpha\succ_k 0$ be concave moduli of continuity. Then every homomorphism $\DD(M;k,\mu)\to\DD(M;k,\alpha)$ has abelian image. If $k\notin\{1,2\}$ then every homomorphism $\DD(M;k,\mu)\to\DD(M;k,\alpha)$ has trivial image.
Moreover, if $k=1$ suppose that $\mu$ is sup-tame, and if $k=2$ suppose that $\mu$ is sub-tame. Then every homomorphism $\DD(M;k,\mu)\to\DD(M;k,\alpha)$ has trivial image.
\end{thm}

In Theorem~\ref{thm:homomorphism}, note that we make no assumptions about continuity of the homomorphisms in question. We remark that Mann~\cite{Mann2015} proved Theorem~\ref{thm:homomorphism} for H\"older continuity moduli $r>s\geq 2$.
Compare this theorem with a more delicate result of Corollary~\ref{cor:cont-hom}.

The proof of Theorem~\ref{thm:homomorphism} is straightforward from Mather's Theorem and Theorem~\ref{thm:spectrum}. Namely, let $\mu$ and $\alpha$ be given and let $Q$ be the group furnished by Theorem~\ref{thm:spectrum}. Then Theorem~\ref{thm:spectrum} shows that $[Q,Q]$ lies in the kernel of any homomorphism $\DD(M;k,\mu)\to\DD(M;k,\alpha)$, so that no such homomorphism is injective. Since $[\DD(M;k,\mu),\DD(M;k,\mu)]$ is simple, the first conclusion of Theorem~\ref{thm:homomorphism} follows. Under the further assumptions on $\mu$, we have that $\DD(M;k,\mu)$ is simple, whence the remaining conclusions of Theorem~\ref{thm:homomorphism} follow.

\subsection{Some open problems}\label{ss:open}

We close by drawing the reader's attention to some open problems which persist in the theory of critical regularity of group actions and in
the context of Mather's Theorem.

\begin{que}\label{que:lip}
Let $k\in\bN$ and let $G$ be a countable subgroup of $\Diff^{k,\mathrm{Lip}}(M)$. Does $G$ admit an injective homomorphism into $\Diff^{k+1}(M)$?
\end{que}

We note that the methods of~\cite{KK2017crit} are not suitable for resolving Question~\ref{que:lip}.
Thurston stability can be used to show that the answer is no for $k=0$~\cite{Calegari2006TAMS}.

\begin{rem}
At the time of the writing of this survey, recent joint work of Kim--Koberda--Mann--Wolff has shown that the answer to Question~\ref{que:lip}
is indeed no for all $k$.
\end{rem}

\begin{que}\label{que:infty}
Is there a countable group $G$ such that $\CR(G)=\infty$ but such that $G$ is not a subgroup of $\Diff^{\infty}(M)$?
\end{que}

It seems likely that the answer to Question~\ref{que:infty} is yes.

\begin{que}\label{que:artin}
Let $A(\gam)$ be a right-angled Artin group of finite critical regularity. What is $\CR(A(\gam))$?
\end{que}

It seems likely that $\CR(A(\gam))$ depends on the structure of $\gam$. Perhaps the simplest right-angled Artin group whose critical regularity is unknown is $(F_2\times\bZ)*\bZ$.

\begin{que}\label{que:fp}
For which $r\in[1,\infty)$ is there a finitely presented group of critical regularity exactly $r$?
\end{que}

We will not speculate as to the answer to Question~\ref{que:fp}. The possible values of $r$ will form a countable set, and $\{1,1.5,2\}$ are all possible values of $r$.

Recall we have proved that $\Diff_+^{1,\alpha}(S^1)$ is perfect whenever $\alpha$ is sup-tame. The following question is not covered by the method described here.
\begin{que}
Is the group $\Diff_+^{1,\alpha}(S^1)$ perfect for $\alpha(s) = s\log(1/s)$?
\end{que}
Note that it is an open question whether or not there exists a non-minimal $C^{1,\alpha}$ action on $S^1$ with an irrational rotation number, for the same $\alpha$~\cite{Herman1979,Navas2011}.
Finally, we reiterate the primary unsolved problem around Mather's Theorem:

\begin{que}
Is the group $\Diff_c^{2}(\bR)$ perfect?
\end{que}

\section*{Acknowledgements}
The first named author is grateful to the organizers of the semester program \emph{Geometry, Topology, and Group Theory in Low Dimensions}, which took place at CIRM (Luminy, France) in 2018. A part of this work was done during his stay for the program. The first and the second named authors are supported by the Mid-Career Researcher Program (2018R1A2B6004003) through the National Research Foundation funded by the government of Korea. The third named author is partially supported by Simons Foundation Collaboration Grant number 429836, by an Alfred P. Sloan Foundation Research Fellowship, and by NSF Grant DMS-1711488. 

\appendix


\providecommand{\bysame}{\leavevmode\hbox to3em{\hrulefill}\thinspace}
\providecommand{\MR}{\relax\ifhmode\unskip\space\fi MR }
\providecommand{\MRhref}[2]{%
  \href{http://www.ams.org/mathscinet-getitem?mr=#1}{#2}
}
\providecommand{\href}[2]{#2}

\end{document}